\documentclass[11pt]{article}

\usepackage{graphicx}
\usepackage{amssymb,amsthm,longtable}
\usepackage{color}
\theoremstyle{definition}

\newtheorem{remark}{Remark}
\theoremstyle{plain}

\newtheorem{proposition}{Proposition}

\newcommand{\ZZ}{\mathbb{Z}}
\newcommand{\NN}{\mathbb{N}}

\newcommand{\ov}{\overline}

\title{Classification of 64-element finite semifields}
\author{ I.F. R\'
ua\thanks{Departamento de Matem\'aticas, Universidad de Oviedo, rua@uniovi.es .
Partially supported by MEC - MTM - 2007 - 67884 Ð C04 - 01}\and El\'\i as F. Combarro\thanks{Artificial Ingelligence Center,
University of Oviedo, \{elias,ranilla\}@aic.uniovi.es . Partially supported by MEC -TIN - 2007 - 61273 and MEC- TIN -2007 - 29664 - E}\and {J. Ranilla$^\dag$}}
\date{}

\begin{document}

\maketitle

\begin{abstract}
A finite semifield $D$ is a finite nonassociative ring with identity such that the set $D^*=D\setminus\{0\}$ is closed under the product. In this paper we obtain a computer-assisted description of all 64-element finite semifields, which completes the classification of finite semifields of order 125 or less.
   
\end{abstract}

\section{Introduction}

A \textbf{finite semifield} (or finite division ring) $D$ is a
finite nonassociative ring with identity such that the set
$D^*=D\setminus\{0\}$ is closed under the product, i.e., it is a loop \cite{Knu65,Cor99}.
Finite semifields have been traditionally considered in the context of finite
geometries since they coordinatize projective semifield planes \cite{Hall}. Recent
applications to coding theory \cite{Calderbank3,codigosdesemicuerpos,Symplectic}, combinatorics and graph theory \cite{ISSAC}, have
broaden the potential interest in these rings.

Because of their diversity, the obtention of general theoretical
algebraic results seems to be a rather difficult (and challenging)
task. On the other hand, because of their finiteness, computational
methods can be naturally considered in the study of these objects. 
So, the classification of finite semifields of a given order is a rather natural problem to use computations. For instance, computers were used in the classification of finite semifields or order 32 \cite{Walker,Knu65} and 81 \cite{Dempwolff}.

In this paper we present a classification up to isotopy of finite semifields with \textbf{64 elements}. Because of the complexity of the problem, the algorithms used in the papers mentioned above can not be directly used to solve the 64-element case. 
 Our techniques combine several old methods and some crucial observations which allow us to reduce dramatically the computational effort. The resulting classification is very rich. Only one tenth of the existing planes was previously known.

The structure of the paper is as follows. In $\S 2$, basic
properties of finite semifields are reviewed. Section $3$ is devoted to known constructions of semifields of order 64.
In $\S
4$, we describe the method we used to classify all 64-element
finite semifields. Finally, in $\S 5$, a complete description of the semifields is given, together with several of their properties.

\section{Preliminaries}

In this section we collect definitions and facts on finite
semifields. Proofs of these results can be found, for instance, in
\cite{Knu65,Cor99}.

    A finite nonassociative ring $D$ is called \textbf{presemifield},
    if the set of nonzero elements $D^*$ is closed under the
    product.
    If $D$ has an identity element, then it is called \textbf{finite semifield}.
If $D$ is a finite semifield, then $D^*$ is a multiplicative loop.
That is, there exists an element $e\in D^*$ (the identity of $D$)
such that $ex=xe=x$, for all $x\in D$ and, for all $a,b\in D^*$,
  the equation $ax=b$
(resp. $xa=b$) has a unique solution.

Apart from finite fields (which are obviously finite semifields), \emph{proper} finite semifields were first considered by L.E. Dickson
\cite{Dickson1} and were deeply studied by A.A. Albert
\cite{Alb52,Alb58,Alb60,Alb61}. The term \emph{finite semifield} was
introduced in 1965 by D.E. Knuth \cite{Knu65}. These rings play an
important role in the study of certain projective planes, called \emph{semifield planes}
\cite{Hall,Knu65}. 

    The characteristic of a finite presemifield $D$ is a prime number $p$, and $D$ is a finite-dimensional algebra over $GF(q)$ ($q=p^c$) of dimension $d$, for some $c,d\in \NN$, so that 
    $|D|=q^d$.
    If $D$ is a finite semifield, then $GF(q)$ can be chosen to be its associative-commutative center $Z(D)$. Other relevant subsets of a finite semifield are the left, right, and middle nuclei ($N_l,N_r,N_m$), and the nucleus $N$.

The definition of isomorphism of  presemifields is the usual one for algebras, and the classification of finite semifields up to isomorphism can be naturally considered. Because of the connections to finite geometries, we must also consider the following notion.
	If $D_1,D_2$ are two presemifields over the same prime field $GF(p)$, then an \textbf{isotopy} between $D_1$ and $D_2$ is a triple $(F,G,H)$ of bijective linear maps $D_1\to D_2$ over $GF(p)$ such that
	$$H(ab)=F(a)G(b)\; \forall a,b\in D_1.$$

It is clear that any isomorphism between two presemifields is an isotopy, but the converse is not necessarily true. Any presemifield is isotopic to a finite semifield \cite[Theorem 4.5.4]{Knu65}.
From any presemifield $D$ a projective plane $\mathcal P(D)$ can be constructed (see \cite{Hall,Knu65} for the details of this construction). Theorem 6 in \cite{Alb60} shows that isotopy of finite semifields is the algebraic translation of the isomorphism between the corresponding projective planes. So, two finite semifields $D_1,D_2$ are isotopic if, and only if, the projective planes $\mathcal P(D_1),\mathcal P(D_2)$ are isomorphic. 

The set of isotopies between a finite semifield $D$ and itself is a group under composition, called the \textbf{autotopy} group, and denoted $\hbox{At}(D)$. This group acts on  the {\it fundamental triangle} of the plane $\mathcal P(D)$, that is, it leaves invariant each of the three lines
$L_x=\{(1,x,0)\ |\ x\in D\}\cup\{(0,1,0)\}\ ,\ L_y=\{(1,0,y)\ |\ y\in D\}\cup\{(0,0,1)\}\ ,\ L_\infty=\{(0,1,z)\ |\ z\in D\}\cup\{(0,0,1)\}.$
	If $D$ is a finite semifield, and $\mathcal D$ is the set of all nonisomorphic semifields isotopic to $D$, then
	$$(|D|-1)^2=|\hbox{At}(D)|\sum_{E\in \mathcal D}\frac{1}{|\hbox{Aut}(E)|}$$
	where $\hbox{Aut}(E)$ denotes the automorphism group of a finite semifield $E$.
	The sum of the right term will be called \emph{the Semifield/Automorphism (S/A) sum} \cite[Theorem 3.3.4]{Knu65}.

If $\mathcal B=[x_1,\dots,x_d]$ is a $GF(q)$-basis of a presemifield $D$, then there exists a unique set of constants $\mathbf A_{D,\mathcal B}=\{A_{i_1i_2i_3}\}_{i_1,i_2,i_3=1}^d\subseteq GF(q)$ such that 
$$x_{i_1}x_{i_2}=\sum_{i_3=1}^d{A_{i_1i_2i_3}}x_{i_3}\; \forall i_1,i_2\in\{1,\dots,d\}$$
This set is called \textbf{3-cube} corresponding to $D$ with respect to the basis $\mathcal B$, and it completely determines the multiplication in $D$.

A remarkable fact is that permutation of the indexes of a 3-cube preserves the absence of nonzero divisors. Namely, if $D$ is a presemifield, and $\sigma\in S_3$ (the symmetric group on the set $\{1,2,3\}$), then the set
$$\mathbf A_{D,\mathcal B}^\sigma=\{A_{i_{\sigma(1)}i_{\sigma(2)}i_{\sigma(3)}}\}_{i_1,i_2,i_3=1}^d\subseteq GF(q)$$ is the 3-cube of a $GF(q)$-algebra $D_{\mathcal B}^\sigma$ which has not zero divisors \cite[Theorem 4.3.1]{Knu65}. 
Notice that, in general, different choice of bases $\mathcal B,\mathcal B'$ lead to nonisomorphic presemifields $D_{\mathcal B}^\sigma,D_{\mathcal B'}^\sigma$. However, these presemifields are always isotopic.

The number of projective planes that can be constructed from a given finite semifield $D$ using the transformation of the group $S_3$ is at most six \cite[Theorem 5.2.1]{Knu65}. Actually, $S_3$ acts on the set of semifield planes of a given order. So, the classification of finite semifields can be reduced to the classification of the corresponding projective planes up to the action of the group $S_3$. In this setting, we will consider a plane as \emph{new} when no \emph{known}\footnote{We have
considered as \emph{known} semifields those appearing in the, up to our knowledge, last survey on the topic, \cite{Kantor}} finite semifield coordinatizes a plane in its $S_3$-class.

 We shall use a graphical representation to distinguish between the different cases. 
The vertices of an hexagon will depict the six different planes obtained from a given finite semifield (cf. \cite[Theorem 5.2.1]{Knu65}).

$$\hbox{ }\hspace{.5cm}\includegraphics[width=3cm,height=3cm]{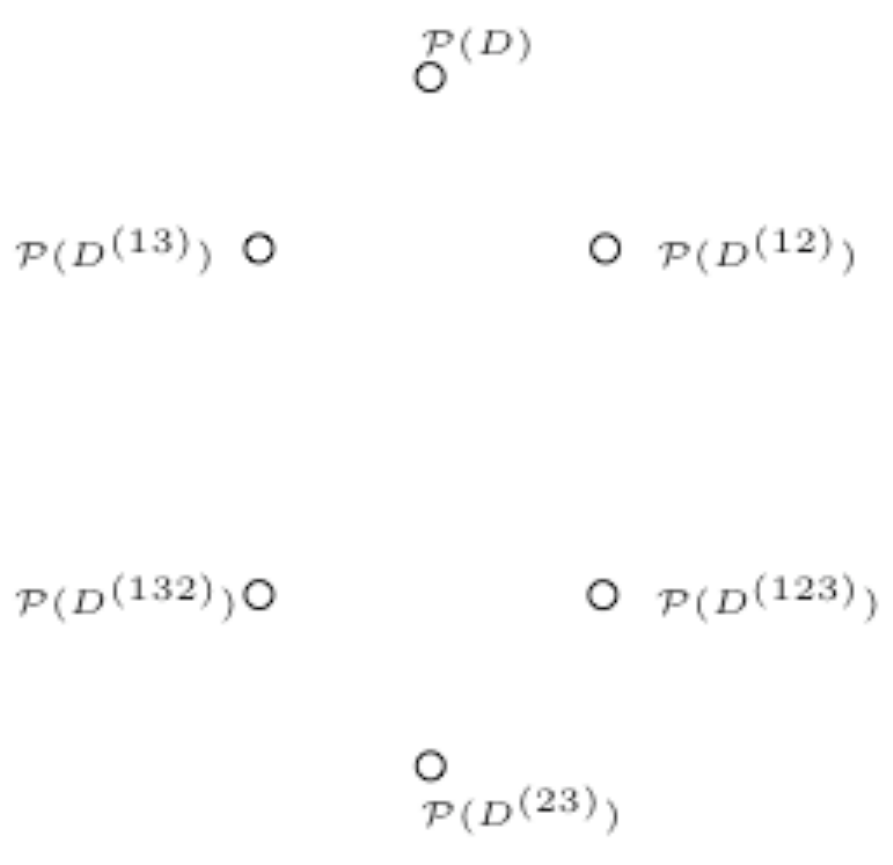}$$

A dotted line between two planes shows that the corresponding finite semifields are isotopic. A continuos line represents that a commutative or symplectic \cite{codigosdesemicuerpos} coordinatizing finite semifield exists.

The construction of finite semifields of a given order can be rephrased
as a matrix problem \cite[Proposition 3]{Hentzel}. We state this proposition in the particular case of semifields of order 64.

\begin{proposition}\label{matrices}
    There exists a finite semifield $D$ of
    64 elements if, and only if, there exists a
    set of $6$ matrices $\{A_1,\dots,A_6\}\subseteq GL(6,2)$ (a \emph{standard basis} \cite{Dempwolff}) such
    that:
    \begin{enumerate}
        \item $A_1$ is the identity matrix;
        \item $\sum_{i=1}^6\lambda_i A_i\in GL(6,2)$, for all
        \textbf{non-zero tuples}
        $(\lambda_1,\dots,\lambda_6)\in  \ZZ_2^d$.
        \item The first column of the matrix $A_i$ is the column vector
         with a $1$ in
    the $i$-th position, and $0$ everywhere else.
    \end{enumerate}
\end{proposition}

This result will be used through this paper to represent 64-element finite semifields. Namely, a finite semifield will be given as a tuple of matrices $(A_2,A_3,A_4,A_5,A_6)$. If the last five columns of
the matrix $A_i$ which has a one in the $i$-th position of the first
column, and zeroes everywhere else are

$$\left(%
\begin{array}{ccccc}
  a_{29} & a_{23} & a_{17} & a_{11} & a_5 \\
  a_{28} & a_{22} & a_{16} & a_{10} & a_4 \\
  a_{27} & a_{21} & a_{15} & a_9 & a_3 \\
  a_{26} & a_{20} & a_{14} & a_8 & a_2 \\
  a_{25} & a_{19} & a_{13} & a_7 & a_1 \\
  a_{24} & a_{18} & a_{12} & a_6 & a_0 \\
\end{array}%
\right)$$
then it is encoded
$\hbox{as the number} \sum_{j=0}^{29}a_j2^j.$

\section{Known semifield planes of order 64}

In this section we give a classification of the known semifield planes of order 64. We considered all the constructions in \cite{Kantor} and explored which planes can be coordinatized by those constructions. This yielded to 31 semifield planes, divided into 11  $S_3$-classes, that we list below.
A semifield representative is given for each plane, together with the order of its automorphism group. Of these planes only two are commutative and can be coordinatized by 14 commutative semifields.

\vspace{0.3cm}\hrule\vspace{0.3cm}

\centerline{\textbf{I (Desarguesian plane)}}

Finite field $\hbox{GF}(64)$ (6 automorphisms) 
$$(A_2,A_3,A_4,A_5,A_6)=
(135274593, 67639409, 33954937, 25632381, 566730623)$$

\centerline{\textbf{II (Twisted field plane)}}

Twisted Field (1 automorphism) with parameters:
\begin{itemize}
	\item $j\in \hbox{GF}(64)$ such that $j^6+j+1=0$;
	\item $\alpha\in \hbox{Aut}(\hbox{GF}(64))$ such that $\alpha(x)=x^4$, for all $x\in \hbox{GF}(64)$;
	\item $\beta\in \hbox{Aut}(\hbox{GF}(64))$ such that $\beta(x)=x^{4^2}$, for all $x\in \hbox{GF}(64)$.
\end{itemize}
$$(A_2,A_3,A_4,A_5,A_6)=
(135274593, 225354480, 673682562, 25632381, 199628676)$$

\centerline{\textbf{III}}

Knuth's semifield of type 2 (1 automorphism) with parameters:
\begin{itemize}
	\item $f\in \hbox{GF}(8)$ such that $f^3+f+1=0$;
	\item $g\in \hbox{GF}(8)$ such that $g+1=0$;
	\item $\sigma\in \hbox{Aut}(\hbox{GF}(8))$ such that $\sigma(x)=x^2$, for all $x\in \hbox{GF}(8)$.
\end{itemize}
$$(A_2,A_3,A_4,A_5,A_6)=
(135274596, 27112887, 35119969, 253266042, 1070246993)$$

\centerline{\textbf{IV}}

Knuth's semifield of type 5 (1 automorphism) with parameters:
\begin{itemize}
	\item $f\in \hbox{GF}(8)$ such that $f^3+f+1=0$;
	\item $g\in \hbox{GF}(8)$ such that $g+1=0$;
	\item $\sigma\in \hbox{Aut}(\hbox{GF}(8))$ such that $\sigma(x)=x^2$, for all $x\in \hbox{GF}(8)$.
\end{itemize}
$$(A_2,A_3,A_4,A_5,A_6)=
(135274593, 189853287, 236639294, 212321269, 624416899)$$

\centerline{\textbf{V (II Huang\&Johnson plane \cite{Huang})} }

Huang\&Johnson sporadic semifield of type II (3 automorphisms).
$$(A_2,A_3,A_4,A_5,A_6)=
(135274623, 1022013944, 102205750, 429859362, 652592216)$$

\centerline{\textbf{VI (III Huang\&Johnson plane)}}

Huang\&Johnson sporadic semifield of type III (3 automorphisms).
$$(A_2,A_3,A_4,A_5,A_6)=
(135274605, 1022014833, 374827988, 557069354, 336124018)$$

\centerline{\textbf{VII (IV Huang\&Johnson plane)}}

Huang\&Johnson sporadic semifield of type IV (1 automorphism).
$$(A_2,A_3,A_4,A_5,A_6)=
(135274605, 427572234, 1072787891, 401402255, 192290736)$$

\centerline{\textbf{VIII (VI Huang\&Johnson plane)}}

Huang\&Johnson sporadic semifield of type VI (2 automorphisms).
$$(A_2,A_3,A_4,A_5,A_6)=
(135274593, 189853287, 580915984, 793113293, 782199145)$$

\centerline{\textbf{IX (VII Huang\&Johnson plane)}}

Huang\&Johnson sporadic semifield of type VII (1 automorphism).
$$(A_2,A_3,A_4,A_5,A_6)=
(135274605, 67640187, 851743451, 194887306, 617256025)$$

\centerline{\textbf{X (VIII Huang\&Johnson plane)}}

Huang\&Johnson sporadic semifield of type VIII (1 automorphism).
$$(A_2,A_3,A_4,A_5,A_6)=
(135274593, 189853287, 1000703633, 930902659, 782199145)$$

\centerline{\textbf{XI}} \centerline{\textbf{(Commutative plane associated to a Kantor-Williams symplectic presemifield)}}

Commutative semifield (6 automorphisms) with tuple of matrices 
$$(A_2,A_3,A_4,A_5,A_6)=
(135274594, 70580276, 37685996, 25345988, 584237329)$$

\vspace{0.3cm}\hrule\vspace{0.3cm}

Apart from the constructions in \cite{Kantor}, the 36 nonprimitive \cite{Wene91} finite semifields of \cite{Hentzel} were also considerd and two new classes (4 planes) were found.

\vspace{0.3cm}\hrule\vspace{0.3cm}

\centerline{\textbf{XII (Two-sided nonprimitive plane)}}

$H$ (Semifield \# 1 in \cite{Hentzel}[page 1423]), the unique nonprimitive semifield of order 64 \cite{Hentzel} (6 automorphisms).
$$(A_2,A_3,A_4,A_5,A_6)=
(146808934, 811798971, 308657185, 563815286, 374228233)$$

\centerline{\textbf{XIII (One-sided nonprimitive plane)}}

Semifield \# 2 in \cite{Hentzel}[page 1423], one-sided nonprimitive semifield (1 automorphism) with tuple of matrices 
$$(A_2,A_3,A_4,A_5,A_6)=
(135274600, 518296613, 253216863, 778190320, 47879003)$$

\vspace{0.3cm}\hrule\vspace{0.3cm}

For each of these 13 semifield representatives, we computed the order of the center and nuclei $ZN=(Z,N,N_l,N_m,N_r)$, the list of all principal isotopes, and the order of their isomorphism groups. Some information on the autotopy group was computed  as well as the length of the orbits in the fundamental triangle $(L_x,L_\infty,L_y)$, given in the form $\sum_{i=1}^ra_i[b_i]$, if $a_i$ cycles of length $b_i$ ($i=1,\dots,r$) exist. All these data are collected in Table \ref{conocidos}.

$ $\hspace{3cm}
\begin{table}
\begin{tabular}{|c|c|c|c|c|c|}
  \hline
   \textbf{Plane}  & $\mathbf{S_3-class}$ & \textbf{$|$At$|$} & $\mathbf{(L_x,L_\infty,L_y)}$&\textbf{S/A sum}& $\mathbf{ZN}$\\
  \hline
  \textbf{I} 	&$\begin{array}{c}\includegraphics[width=1cm,height=1cm]{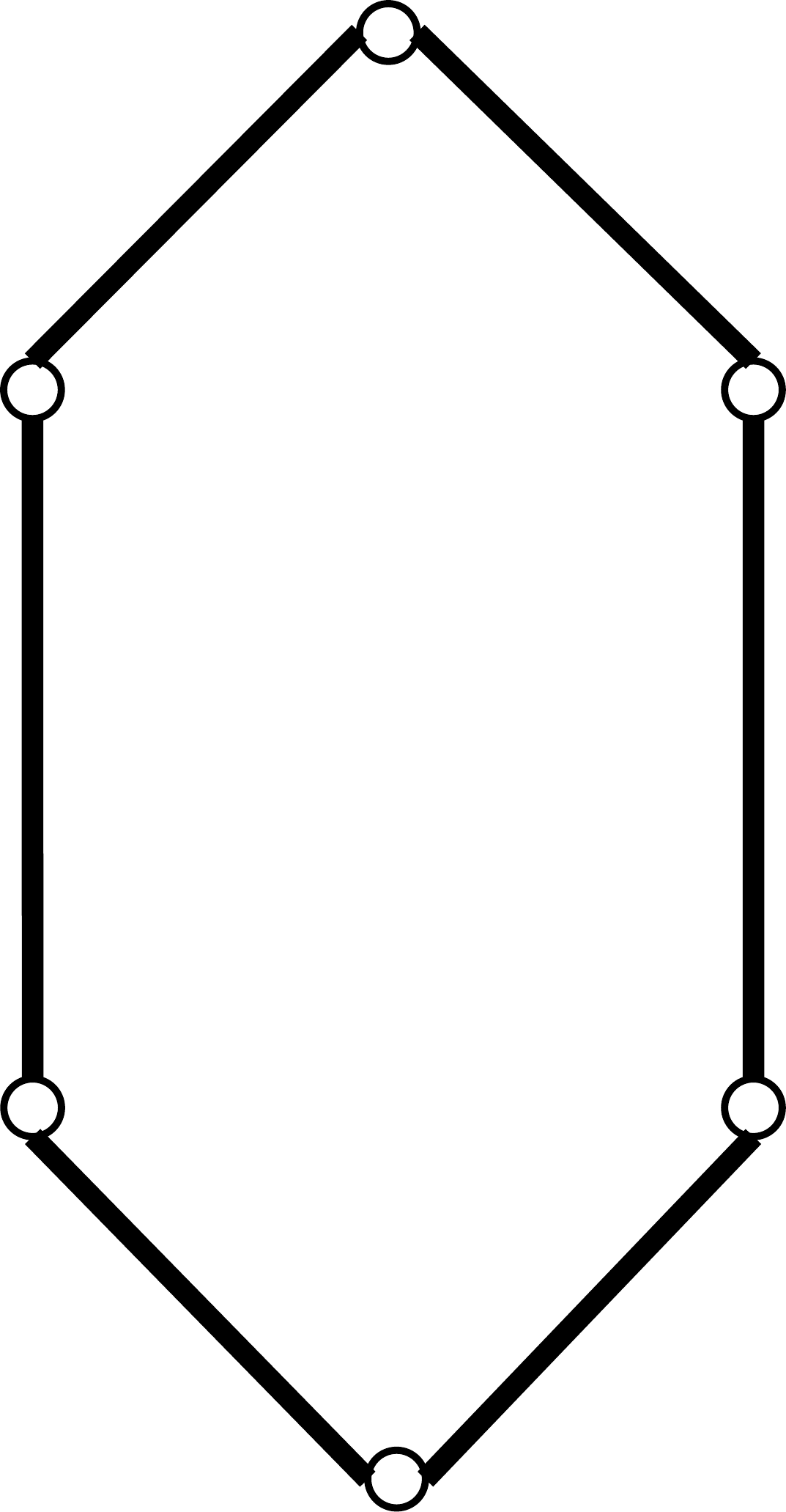}\\\end{array}$&	
  $\begin{array}{c}    23814\\   
  % \hbox{Nilp. exp. 4}
  \\\end{array}$ 
  & $\begin{array}{c}    2[1]+1[80]\\    2[1]+1[80]\\    2[1]+1[80]\\\end{array}$  	 & $\frac{1}{6}$	& $(64,64,64,64,64)$\\\hline
\textbf{II} 	&$\begin{array}{c}\includegraphics[width=1cm,height=1cm]{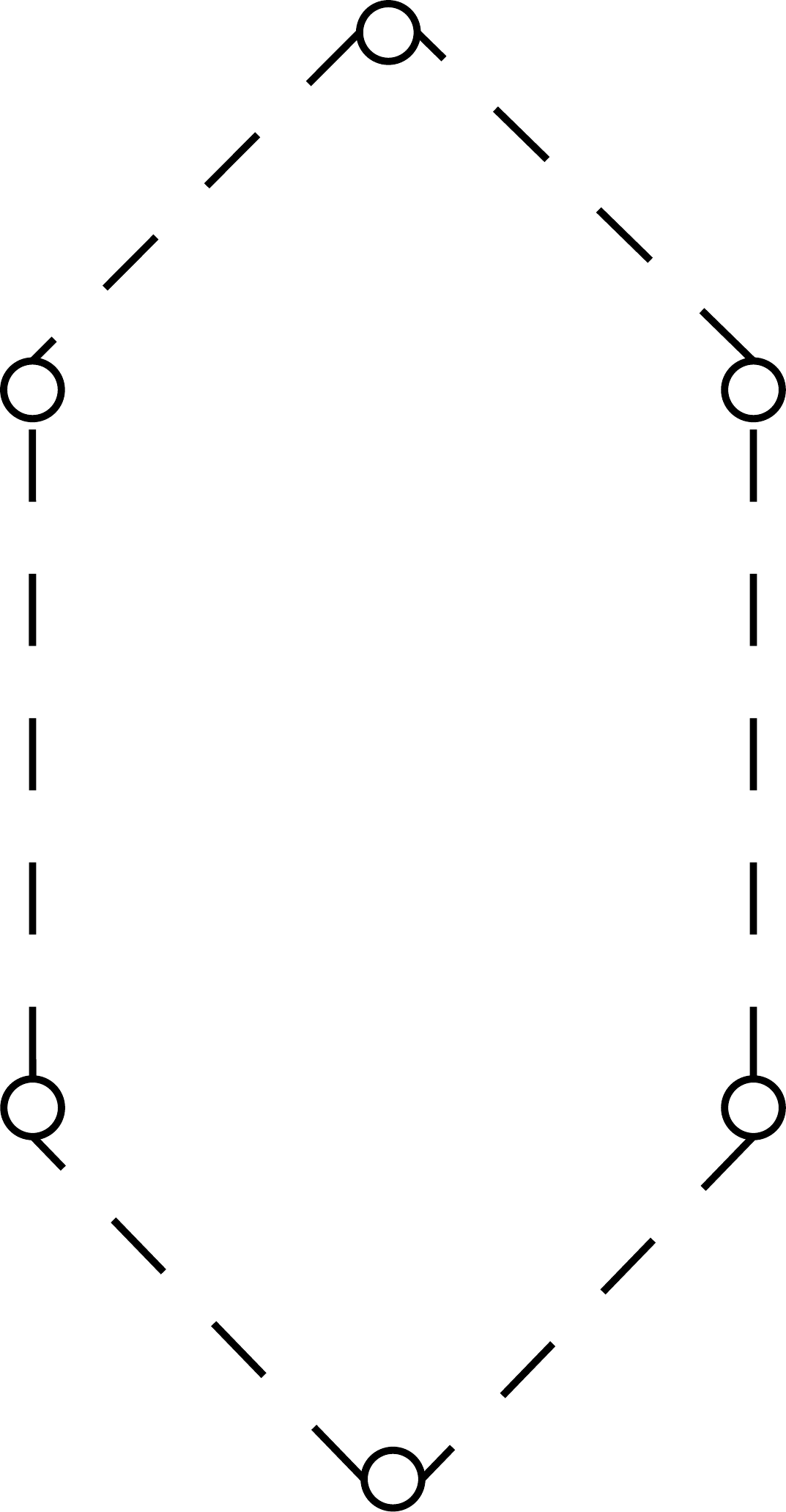}\\\end{array}$&	$  \begin{array}{c} 567 \\ \hbox{Solvable}\\\end{array}$ & $\begin{array}{c} 2[1]+1[63]\\ 2[1]+1[63]\\ 2[1]+1[63]\\\end{array}$	 & $\frac{7}{1}$	& $(4,4,4,4,4)$\\\hline
\textbf{III} 	&$\begin{array}{c}\includegraphics[width=1cm,height=1cm]{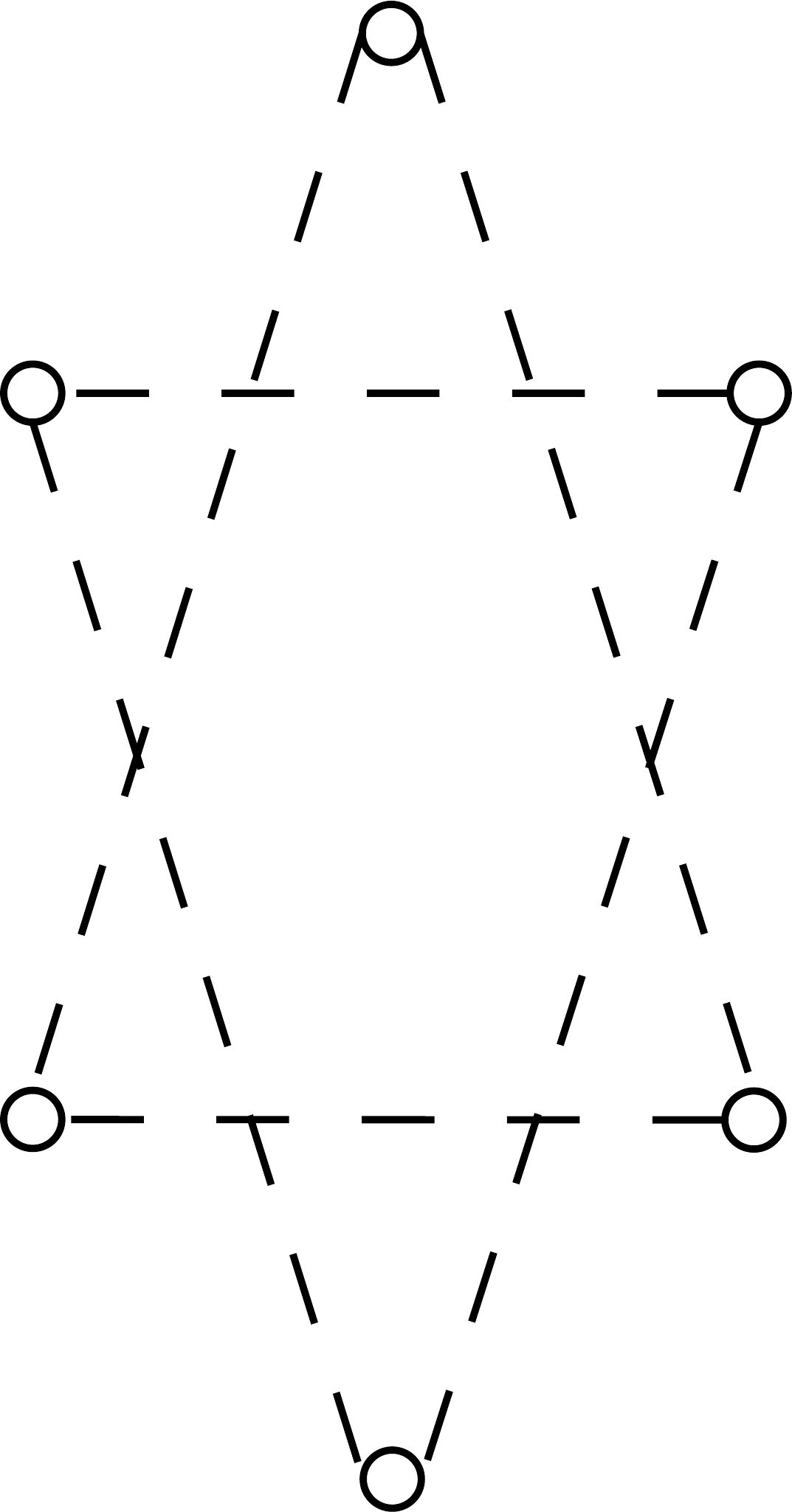}\\\end{array}$&	$  \begin{array}{c} 49 \\ \hbox{$\ZZ_7^2$}\\\end{array}$ & $\begin{array}{c} 2[1]+2[7]+1[49]\\ 2[1]+2[7]+1[49]\\ 2[1]+2[7]+1[49]\\\end{array}$	 & $\frac{81}{1}$	& $(2,2,2,2,2)$\\\hline
\textbf{IV} 	&$\begin{array}{c}\includegraphics[width=1cm,height=1cm]{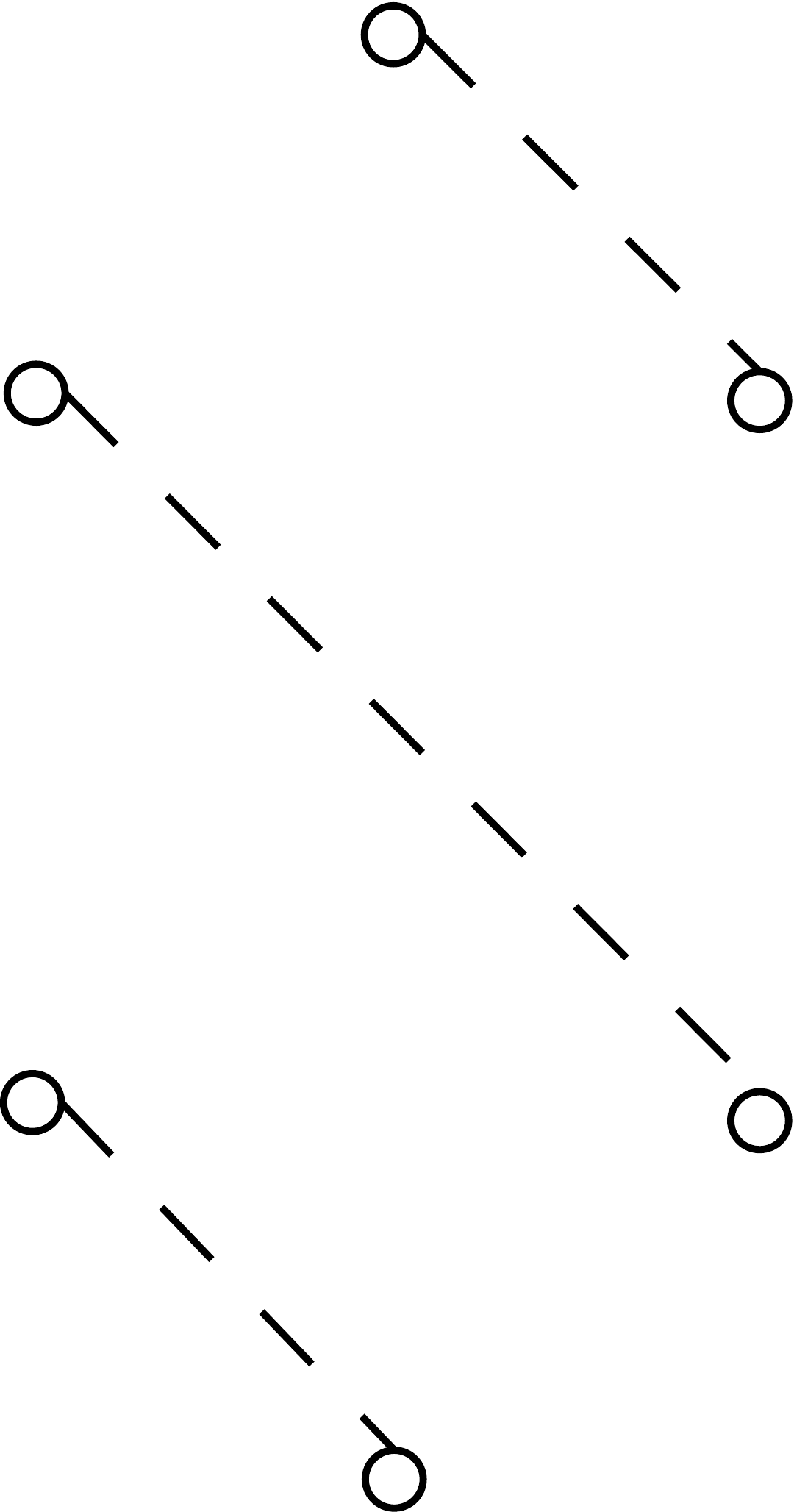}\\\end{array}$&	$  \begin{array}{c} 441 \\ \hbox{Solvable}\\\end{array}$ & $\begin{array}{c} 2[1]+1[63]\\ 2[1]+1[63]\\ 2[1]+2[7]+1[49]\\\end{array}$	 & $\frac{9}{1}$	& $(2,2,8,4,8)$\\\hline
\textbf{V} 	&$\begin{array}{c}\includegraphics[width=1cm,height=1cm]{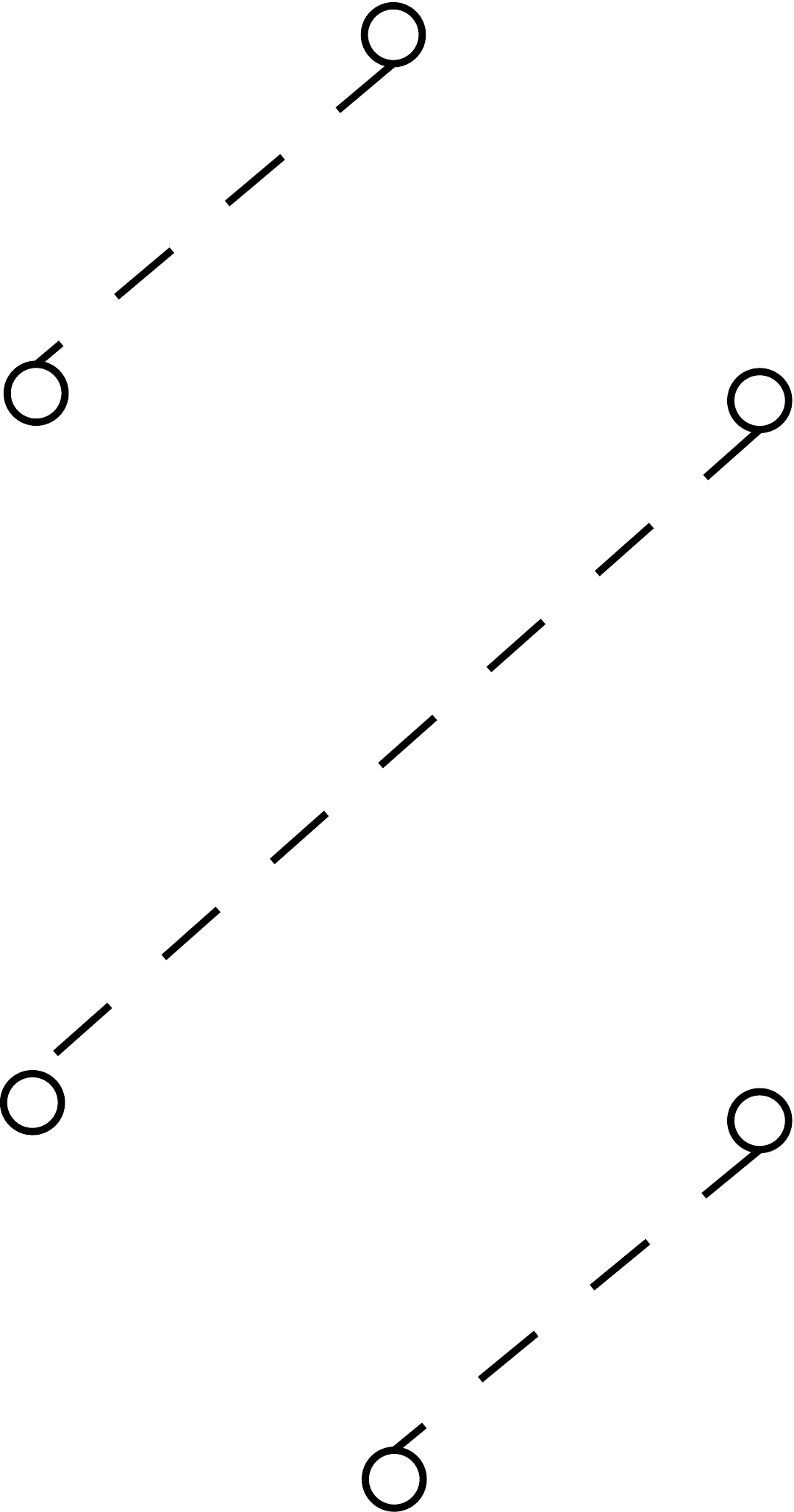}\\\end{array}$&	$  \begin{array}{c} 42 \\ \hbox{Solvable}\\\end{array}$ & $\begin{array}{c} 3[1]+1[2]+2[3]+9[6]\\ 2[1]+1[7]+1[14]+1[42]\\ 2[1]+1[7]+1[14]+1[42]\\\end{array}$	 & $\frac{92}{1} + \frac{2}{2} + \frac{4}{3} + \frac{1}{6}$	& $(2,2,8,2,2)$\\\hline
\textbf{VI} 	&$\begin{array}{c}\includegraphics[width=1cm,height=1cm]{hexagonLinTrB.pdf}\\\end{array}$&	$  \begin{array}{c} 42 \\ \hbox{Solvable}\\\end{array}$ & $\begin{array}{c} 3[1]+1[2]+2[3]+9[6]\\ 2[1]+1[7]+1[14]+1[42]\\ 2[1]+1[7]+1[14]+1[42]\\\end{array}$	 & $\frac{92}{1} + \frac{2}{2} + \frac{4}{3} + \frac{1}{6}$	& $(2,2,8,2,2)$\\\hline
\textbf{VII} 	&$\begin{array}{c}\includegraphics[width=1cm,height=1cm]{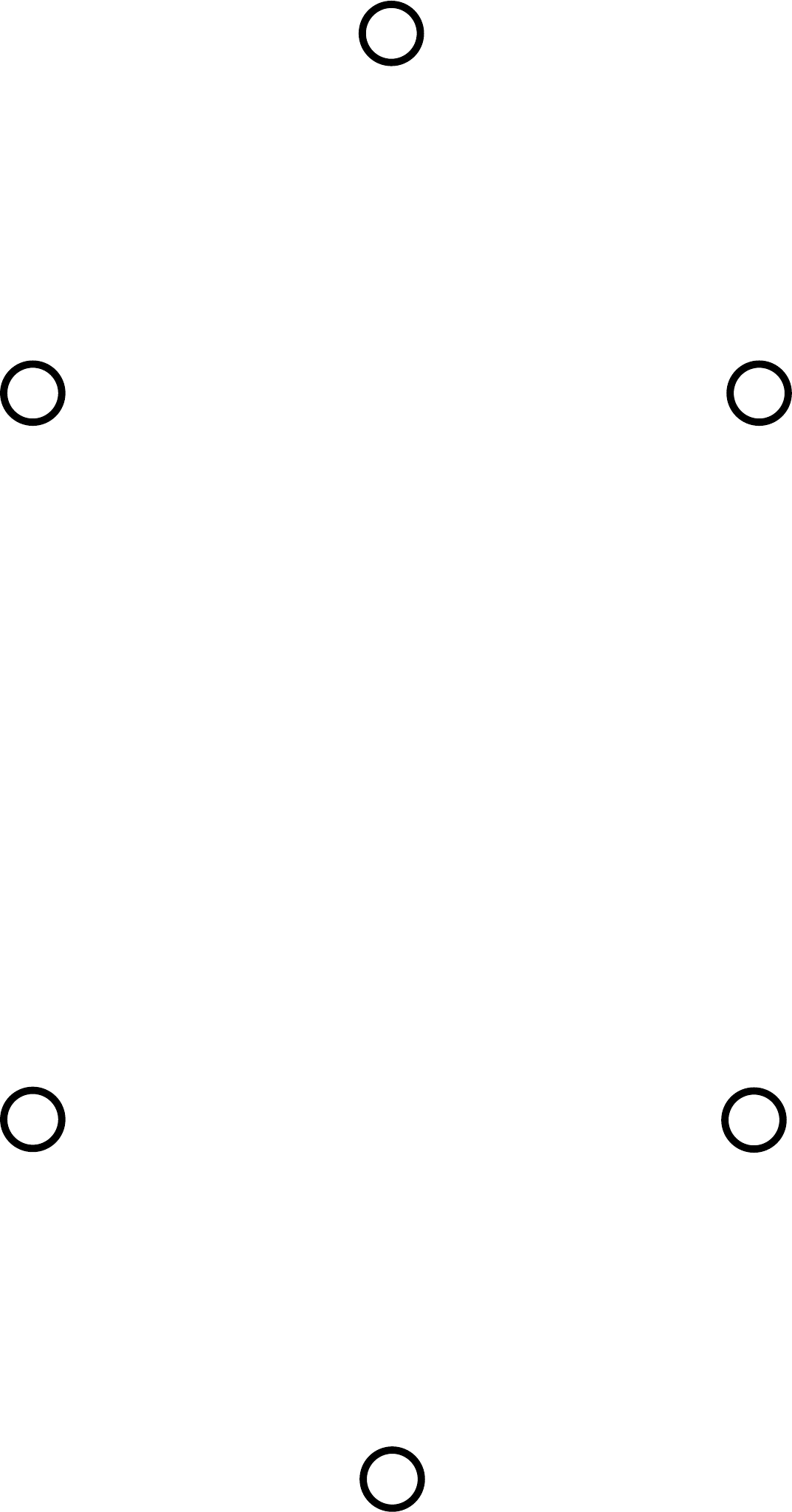}\\\end{array}$&	$  \begin{array}{c} 14 \\ \hbox{Cyclic}\\\end{array}$ & $\begin{array}{c} 9[1]+28[2]\\ 2[1]+1[7]+4[14]\\ 2[1]+1[7]+4[14]\\\end{array}$	 & $\frac{280}{1} + \frac{7}{2}$	& $(2,2,8,2,2)$\\\hline
\textbf{VIII} 	&$\begin{array}{c}\includegraphics[width=1cm,height=1cm]{hexagonLinTrB.pdf}\\\end{array}$&	$  \begin{array}{c} 126 \\ \hbox{Solvable}\\\end{array}$ & $\begin{array}{c} 2[1]+4[3]+1[6]+3[9]+1[18]\\ 2[1]+1[21]+1[42]\\ 2[1]+1[21]+1[42]\\\end{array}$	 & $\frac{28}{1} + \frac{7}{2}$	& $(2,2,8,4,4)$\\\hline
\textbf{IX} 	&$\begin{array}{c}\includegraphics[width=1cm,height=1cm]{hexagonLinTrB.pdf}\\\end{array}$&	$  \begin{array}{c} 42 \\ \hbox{Solvable}\\\end{array}$ & $\begin{array}{c} 3[1]+1[2]+6[3]+7[6]\\ 2[1]+1[21]+1[42]\\ 2[1]+1[21]+1[42]\\\end{array}$	 & $\frac{91}{1} + \frac{7}{2}$	& $(2,2,8,2,2)$\\\hline
\textbf{X} 	&$\begin{array}{c}\includegraphics[width=1cm,height=1cm]{hexagonLinTrB.pdf}\\\end{array}$&	$  \begin{array}{c} 42 \\ \hbox{Solvable}\\\end{array}$ & $\begin{array}{c} 5[1]+6[2]+4[3]+6[6]\\ 2[1]+1[21]+1[42]\\ 2[1]+1[21]+1[42]\\\end{array}$	 & $\frac{91}{1} + \frac{7}{2}$	& $(2,2,8,2,2)$\\\hline
\textbf{XI} 	&$\begin{array}{c}\includegraphics[width=1cm,height=1cm]{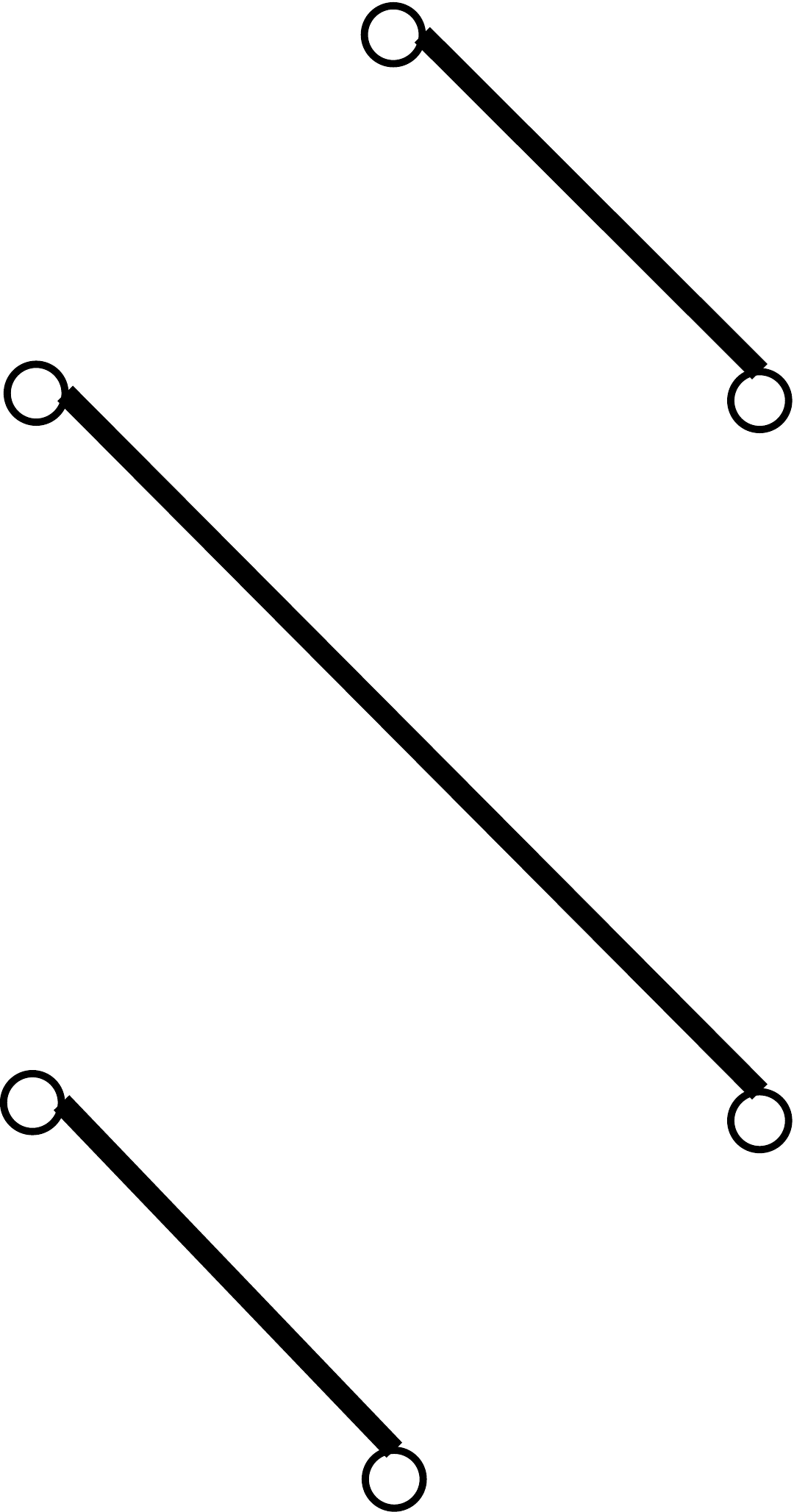}\\\end{array}$&	$  \begin{array}{c} 18 \\ \hbox{Solvable}\\\end{array}$ & $\begin{array}{c} 2[1]+1[3]+1[6]+2[9]+2[18]\\ 2[1]+1[3]+1[6]+2[9]+2[18]\\ 3[1]+1[2]+2[3]+9[6]\\\end{array}$	 & $\frac{211}{1} + \frac{16}{2} + \frac{4}{3} + \frac{1}{6}$	& $(2,2,2,4,2)$\\\hline
\textbf{XII} 	&$\begin{array}{c}\includegraphics[width=1cm,height=1cm]{hexagonCompTr.pdf}\\\end{array}$&	$  \begin{array}{c} 6 \\ \hbox{Symmetric}\\\end{array}$ & $\begin{array}{c} 3[1]+1[2]+6[3]+7[6]\\ 3[1]+1[2]+6[3]+7[6]\\ 3[1]+1[2]+6[3]+7[6]\\\end{array}$	 & $\frac{636}{1} + \frac{48}{2} + \frac{4}{3} + \frac{1}{6}$	& $(2,2,2,2,2)$\\\hline
\textbf{XIII} 	&$\begin{array}{c}\includegraphics[width=1cm,height=1cm]{hexagonLinTrB.pdf}\\\end{array}$&	$  \begin{array}{c} 3 \\ \hbox{Cyclic}\\\end{array}$ & $\begin{array}{c} 17[1]+16[3]\\ 2[1]+21[3]\\ 2[1]+21[3]\\\end{array}$	 & $\frac{1323}{1}$	& $(2,2,2,2,2)$\\\hline

\end{tabular}
\caption{Known finite semifields of 64 elements and their properties}\label{conocidos}\end{table}%

\begin{remark}\label{primitividad} 
\begin{enumerate}
\item Let us notice that the  Huang\&Johnson plane of type V \cite{Huang} is $S_3$-equivalent to the plane \textbf{VII} above. Namely, it is the plane \textbf{VII}$^{(23)}$. 
\item
All these planes can be coordinatized by a finite semifield containing a primitive element (even the planes \textbf{XII} and \textbf{XIII}).
\end{enumerate}
\end{remark}

Let us introduce presentations of the autotopism group of some semifields in the table above. 
$$\hbox{At}_{\hbox{\textbf{II}}}=<x,y,z,t,u \ |\ x^3 = z, 
    y^3 =t^2, 
    z^3 =  
    t^3 =  
    u^7 = 1, 
    y^x = yt, 
    u^x = u^2>$$
$$\hbox{At}_{\hbox{\textbf{IV}}}=<x,y,z,t\ |\ 
x^3 = y, 
    y^3 =  
    z^7 = 
    t^7 = 1, 
    z^x = z^2, 
    t^x = t^2>$$

$$\hbox{At}_{\hbox{\textbf{V}}}=\hbox{At}_{\hbox{\textbf{VI}}}=<x,y,z\ |\ 
x^2 =  
    y^3 =  
    z^7 = 1, 
    z^y = z^2>$$

$$\hbox{At}_{\hbox{\textbf{VIII}}}=<x,y,z,t\ |\ 
x^2 =  
    y^7 =  
    z^3 = 
    t^3 = 1, 
    z^x = z^2, 
    t^x = t^2>$$
    
    $$\hbox{At}_{\hbox{\textbf{IX}}}=\hbox{At}_{\hbox{\textbf{X}}}=<x,y,z\ |\ 
x^2 =  
    z^3 =  
    y^7 = 1, 
    z^x = z^2>$$
    
    $$\hbox{At}_{\hbox{\textbf{XI}}}=<x,y,z\ |\ 
x^2 =  
    y^3 =  
    z^3 = 1, 
    z^x = z^2>$$

Since the same plane can be coordinatized by several constructions, we include Table \ref{coordenadores}, where every plane is followed by a list of those (pre)semifield constructions that coordinatize it:

\begin{itemize}
\item FF: Finite Field
\item TF: Twisted Field
\item K: Knuth's semifield of types 1 to 5 
\item JJ: Jha-Johnson semifield constructed over $GF(4)$ or $GF(8)$
\item HJ: Huang-Johnson esporadic semifield II to VII
\item KW: Kantor-Williams symplectic presemifield or associated commutative presemifield
\item NP: Nonprimitive semifield (1 or 2-sided)
\end{itemize}

\begin{table}[htdp]
\begin{center}
\begin{tabular}{|c||c|c|c|c|c|c|}
\hline
Plane $\backslash$ Permutation& 1 & (1 2) & (1 3) & (1 2 3) & (2 3) & (1 3 2)\\
\hline
\hline
I& FF& --&--&--&--&--\\
\hline
II&TF&--&--&--&--&--\\
\hline
III&K2&K2&--&--&--&--\\
\hline
IV&K5&--&K3 / JJ8 &--&K4&--\\
\hline
V&HJII&NONE&--&NONE&--&--\\
\hline
VI&HJIII&NONE&--&NONE&--&--\\
\hline
VII&HJIV&--&--&--&HJV&--\\
\hline
VIII&HJVI / JJ4&NONE&--&NONE&--&--\\
\hline
IX&HJVII&NONE&--&NONE&--&--\\
\hline
X&HJVIII&NONE&--&NONE&--&--\\\hline
XI&KWC&--&KWS&--&NONE&--\\\hline
XII&NP2&--&--&--&--&--\\\hline
XIII&NP1&NP1&--&NONE&--&--\\
\hline
\end{tabular}
\caption{Known planes and their coordinatizing constructions}\label{coordenadores}
\end{center}
\end{table}%

\section{Search algorithm}

In this section we describe an algorithm that generates all finite semifields of order 64. As we have previously observed (Proposition \ref{matrices}), any 64-element semifield can be described as a tuple of 6 matrices satisfying certain conditions. So, the output of the algorithm will be tuples of matrices which correspond to finite semifields. Not all possible tuples satisfying the conditions of Proposition \ref{matrices} will be listed.
It is only necessary to obtain representatives of all $S_3$-classes of equivalence. 

Our method is based on the algorithm proposed in \cite{Hentzel}. This algorithm proved to be extraordinary efficient in the computation of all nonprimitive semifields of order 64, and it can be easily adapted to our case. It first takes matrices $A_1$ and $A_2$, the first equal to the identity matrix, the second one chosen among a small amount of matrices (we will comment on this latter). Then, it computes 15 lists $L_{(\lambda_3,\lambda_4,\lambda_5,\lambda_6)}$ (where $(\lambda_3,\lambda_4,\lambda_5,\lambda_6)\in \ZZ_2^4\setminus\{(0,0,0,0)\}$) containing matrices $B$ of $GL(6,2)$ such that the first column of $B$ is the vector $(0,0,\lambda_3,\lambda_4,\lambda_5,\lambda_6)^t$, and such that the matrices $B+A2, B+A1$ and $B+A2+A1$ are elements of $GL(6,2)$. If $(A_1,A_2,A_3,A_4,A_5,A_6)$ is a tuple corresponding to a finite semifield of 64 elements, then any nonzero linear combination $\sum_{i=3}^6\lambda_iA_i$ must be contained in the corresponding list 
$L_{(\lambda_3,\lambda_4,\lambda_5,\lambda_6)}$. This provides a fast test to check if a tuple of matrices corresponds to a finite semifield. 

Then, these lists are compared to produce consistent tuples of matrices $(A_3,A_4,A_5,A_6)$. This procedure, which is described in detail in \cite{Hentzel}, essentially sieves the lists of matrices $A_4,A_5$ and $A_6$ with the help of the other lists (the matrix $A_3$ is fixed at this point). The resulting tiny lists are used to generate tuples of  6 matrices, and it is checked whether these tuples are related to a finite semifield or not. Once all these tuples are generated, a classification algorithm produces the different isomorphism classes of nonprimitive finite semifields.

This method produced satisfactory results in \cite{Hentzel}, even though the  computational effort was remarkable big (if run in a single machine, 10 months of computing time). However, it can not be directly applied to our case. The main obstruction is the size of the lists $L_{(\lambda_3,\lambda_4,\lambda_5,\lambda_6)}$. The lists created in \cite{Hentzel}, because of some extra conditions, contained approximately 2 million matrices each. In the new situation, the lists contain more than 7 million matrices each. If we applied directly the method, it would take approximately 466 months to complete the task. This means that, in a certain sense, the classification problem of 64-element finite semifields is almost 50 times more difficult than the primitivity problem.
The main drawback of the method is that all tuples of matrices must be computed, and classification is only achieved after search.

The algorithm which was used in \cite{Dempwolff} to classify 81-element finite semifields combines search and classification. It works as follows. Tuples of four matrices $(A_1,A_2,A_3,A_4)$ associated to finite semifields of order 81 are obtained in 4 steps. In the $i$-th step a list of \emph{nonequivalent partial standard bases of size $i$}  is computed \cite{Dempwolff}. This means that a list of tuples $(A_1,\dots,A_i)$ of $i$ matrices is produced, with the following two conditions. First, any tuple in the list can be potentially extended to get a tuple of 4 matrices associated to a finite semifield. Second, none of them can be obtained from any other by means of an isotopy transformation. This ensures that the resulting tuples of 4 matrices coordinatize nonisotopic finite semifields. This method lead to the classification of 81-element finite semifields \emph{in a few days on a PC} \cite{Dempwolff}.

This algorithm can not be directly used in our setting either. First, because it is necessary to obtain tuples of 6 matrices, which represents 6 steps instead of 4. Moreover, these matrices have order $6\times 6$, and not $4\times 4$. This means a significant increase of computation time. Also, some storing problems do appear in the intermediate steps, since the number of nonequivalent partial standard bases of size 4 is huge. 

We managed to combine these two methods in an effective way. Namely, we used first the  method of \cite{Dempwolff} to produce nonequivalente partial standard bases of size 3.  Since $S_3$-classes are enough for the classification of finite semifields, we considered nonequivalent partial standard bases up to $S_3$-equivalence: the matrices of a partial standard basis of a semifield $D$, when transposed, are a partial standard basis of $D^{(1,3)}$. 
Then, the algorithm of \cite{Hentzel} was used to complete these bases to tuples of 6 matrices. 
This combination of methods produced satisfactory results since there was a 95\% reduction in the number of cases to be explored. Specifically, a total amount of 399866 nonequivalent partial standard bases of size 3 were considered. Of these, 377675 had matrix $A_2=C(x^6+x+1)$ (companion matrix), and the rest had matrix $A_2=C(x^6+x^5+x^3+x^2+1)$. Notice that, because of \cite[Section 5]{Hentzel} and Remark \ref{primitividad} above, the matrix $A_2$ can be always chosen to be a companion matrix of a primitive polynomial.
The implemented algorithm in language C required 30 days on a 12 2.5 GHz CPU linux cluster (1 computer year on a single PC). The output consisted in  95877 tuples. Classification of these matrices up to $S_3$-equivalence lead to the results of the next section. 

Let us remark that the extraordinary feature of our algorithm is that it allowed us to solve the classification problem in approximately the same time it was needed in \cite{Hentzel} to solve the primitivity problem. As noticed above, this problem is, in a certain sense, 50 times simpler.

\section{New Semifield Planes of order 64: a classification}

The classification of 64-element finite semifields that we present in this section, completes the classification of finite semifields of order 125 or less \cite{Kleinfeld,Walker,Knu65,Dempwolff,Meni1,Meni2}. Let us compare the number of $S_3$-equivalence classes, semifield planes, and coordinatizing finite semifields which were found, with those previously known (Table \ref{numeros}).

\begin{table}[htdp]
\begin{center}
\begin{tabular}{|c|c|c|c|}
  \hline
   \textbf{Number of classes}  & \textbf{$S_3$-action} & \textbf{Isotopy} & \textbf{Isomorphism}\\
  \hline
  Previously known &  13&35&3058 \\

  \hline
   Actual number & 80&332&87714 \\
   \hline
\end{tabular}
\caption{Number of 64-element finite semifields}\label{numeros}
\end{center}
\end{table}%

As we can see approximately one tenth of the semifield planes which exist were previously known.
The matrices of the coordinatizing semifields, from $S_3$-classes XIV to LXXX are collected in Table \ref{numerajos}. A semifield representative with maximal number of automorphisms was chosen for each $S_3$-class.

\begin{center}
\begin{longtable}{|c|c|c|c|c|c|}
\caption{Matrices of new 64-element finite semifields}\label{numerajos} \\
\hline  $\#$ & $A_2$  & $A_3$ & $A_4$ & $A_5$ & $A_6$\\ 
\endfirsthead

\multicolumn{6}{c}%
{{ \tablename\ \thetable{} -- continued from previous page}} \\\hline  $\#$ & $A_2$  & $A_3$ & $A_4$ & $A_5$ & $A_6$\\ \hline
\endhead

\hline
\endfoot

\hline 
\endlastfoot
 
  \hline\textbf{	XIV	}	& 135274594 & 833399958 & 260289148 & 1031543734 & 289062724 \\
\textbf{	XV	}	& 135274608 & 929017139 & 43661225 & 236902583 & 58939658 \\
\textbf{	XVI	}	& 135274599 & 780494300 & 790514579 & 230544263 & 113930782 \\
\textbf{	XVII	}	& 135274613 & 129115169 & 42444851 & 21901924 & 898304298 \\
\textbf{	XVIII	}	& 135274600 & 61518615 & 173580196 & 594359470 & 490618435 \\
\textbf{	XIX	}	& 135274596 & 127101463 & 659059070 & 333784654 & 916173704 \\
\textbf{	XX	}	& 135274594 & 71778966 & 584275037 & 294214292 & 47830349 \\
\textbf{	XXI	}	& 135274623 & 210604913 & 369031692 & 869572955 & 513934562 \\
\textbf{	XXII	}	& 135274614 & 298719213 & 371084575 & 198518457 & 1064230576 \\
\textbf{	XXIII	}	& 135274593 & 861665485 & 782365624 & 171288969 & 901459391 \\
\textbf{	XXIV	}	& 135274599 & 329921844 & 727137562 & 76464013 & 227638817 \\
\textbf{	XXV	}	& 135274611 & 502051974 & 921041233 & 799525002 & 1031595865 \\
\textbf{	XXVI	}	& 135274608 & 245456463 & 960054086 & 892794578 & 298620733 \\
\textbf{	XXVII	}	& 135274605 & 685174186 & 978351539 & 30452336 & 770639372 \\
\textbf{	XXVIII	}	& 135274605 & 287116498 & 229601764 & 116769706 & 659794001 \\
\textbf{	XXIX	}	& 135274600 & 391353891 & 935896110 & 611263392 & 219544639 \\
\textbf{	XXX	}	& 135274605 & 253357853 & 818841952 & 214536771 & 376932474 \\
\textbf{	XXXI	}	& 135274605 & 927823728 & 1043775209 & 180241271 & 355149199 \\
\textbf{	XXXII	}	& 135274594 & 830091251 & 496454571 & 294214292 & 937736232 \\
\textbf{	XXXIII	}	& 135274594 & 625566993 & 820897994 & 1041833019 & 527914637 \\
\textbf{	XXXIV	}	& 135274605 & 1015682606 & 420972778 & 331973660 & 602075080 \\
\textbf{	XXXV	}	& 135274605 & 760830577 & 942032486 & 331689540 & 1059713288 \\
\textbf{	XXXVI	}	& 135274600 & 1013759534 & 967035803 & 19745382 & 54171530 \\
\textbf{	XXXVII	}	& 135274618 & 702966560 & 581719755 & 551123260 & 248773288 \\
\textbf{	XXXVIII	}	& 135274594 & 615209500 & 653723442 & 597803524 & 64203530 \\
\textbf{	XXXIX	}	& 135274594 & 766078631 & 695326175 & 915348146 & 1069733783 \\
\textbf{	XL	}	& 135274623 & 151599276 & 350755192 & 628485436 & 451306380 \\
\textbf{	XLI	}	& 135274603 & 977768416 & 817386359 & 868142796 & 1019542329 \\
\textbf{	XLII	}	& 135274618 & 815130072 & 525371889 & 729166901 & 111689247 \\
\textbf{	XLIII	}	& 135274606 & 106019670 & 834133431 & 978786631 & 810064469 \\
\textbf{	XLIV	}	& 135274605 & 127050072 & 1067943835 & 187195280 & 46720452 \\
\textbf{	XLV	}	& 135274611 & 1029317650 & 520088048 & 219035419 & 1059375133 \\
\textbf{	XLVI	}	& 135274593 & 331766335 & 841368844 & 1036592040 & 609803946 \\
\textbf{	XLVII	}	& 135274608 & 153781051 & 728397374 & 250468680 & 37621084 \\
\textbf{	XLVIII	}	& 135274613 & 194046320 & 574609178 & 123846514 & 234580732 \\
\textbf{	XLIX	}	& 135274606 & 769446985 & 220872112 & 247312724 & 317866821 \\
\textbf{	L	}	& 135274620 & 658083501 & 394974963 & 168700243 & 63130518 \\
\textbf{	LI	}	& 135274623 & 834590458 & 173628128 & 932232776 & 129367971 \\
\textbf{	LII	}	& 135274623 & 836224690 & 207465030 & 645936619 & 243736264 \\
\textbf{	LIII	}	& 135274613 & 799452523 & 265350121 & 339974318 & 626807932 \\
\textbf{	LIV	}	& 135274617 & 1021351255 & 508837061 & 331738527 & 774541306 \\
\textbf{	LV	}	& 135274603 & 685722191 & 420150503 & 99266464 & 259409008 \\
\textbf{	LVI	}	& 135274594 & 256099125 & 620378737 & 948416852 & 976951193 \\
\textbf{	LVII	}	& 135274596 & 709090705 & 611750851 & 500630203 & 252344113 \\
\textbf{	LVIII	}	& 135274599 & 296969012 & 813505114 & 902870605 & 253337313 \\
\textbf{	LIX	}	& 135274620 & 1001102903 & 554659815 & 887933802 & 199592399 \\
\textbf{	LX	}	& 135274617 & 1033514144 & 124586730 & 857117276 & 1057874229 \\
\textbf{	LXI	}	& 135274596 & 24792404 & 923093719 & 241684897 & 311785609 \\
\textbf{	LXII	}	& 135274600 & 335867443 & 64403431 & 539538757 & 1048273875 \\
\textbf{	LXIII	}	& 135274599 & 379947884 & 468255421 & 166163471 & 737134591 \\
\textbf{	LXIV	}	& 135274608 & 675671366 & 173741043 & 1064349970 & 330692042 \\
\textbf{	LXV	}	& 135274599 & 773727160 & 109764209 & 24488523 & 96297681 \\
\textbf{	LXVI	}	& 135274603 & 116650580 & 995749484 & 402467592 & 539691559 \\
\textbf{	LXVII	}	& 135274617 & 17661426 & 345529328 & 1017107050 & 175067510 \\
\textbf{	LXVIII	}	& 135274618 & 213725930 & 722775913 & 632723229 & 390126804 \\
\textbf{	LXIX	}	& 135274614 & 298218413 & 614428189 & 791447010 & 595936461 \\
\textbf{	LXX	}	& 135274617 & 885643496 & 663907199 & 176085274 & 942663142 \\
\textbf{	LXXI	}	& 135274620 & 508419752 & 1024350995 & 67020677 & 581629534 \\
\textbf{	LXXII	}	& 135274606 & 559395716 & 262131527 & 86612471 & 850591393 \\
\textbf{	LXXIII	}	& 135274593 & 808819530 & 46309136 & 231311151 & 1016044841 \\
\textbf{	LXXIV	}	& 135274608 & 382084651 & 723363276 & 222713684 & 40992056 \\
\textbf{	LXXV	}	& 135274605 & 844321142 & 295097610 & 835339006 & 805165097 \\
\textbf{	LXXVI	}	& 135274617 & 342415979 & 844418594 & 360800935 & 877880127 \\
\textbf{	LXXVII	}	& 135274596 & 61984053 & 93250846 & 172201558 & 898164556 \\
\textbf{	LXXVIII	}	& 135274613 & 844649372 & 514653314 & 956814412 & 890109555 \\
\textbf{	LXXIX	}	& 135274608 & 580417165 & 384216079 & 701045724 & 922685158 \\
\textbf{	LXXX	}	& 135274614 & 599051760 & 670195531 & 755889110 & 1021850782 \\
\end{longtable}
\end{center}

We processed these semifield representatives to obtain the order of the center and nuclei, the list of all principal isotopes, and the order of their isomorphism groups. Also, the length of the orbits in the fundamental triangle and some information on the autotopy group was computed. All these data are collected in Table \ref{superdatos}.

\begin{center}
\begin{longtable}{|c|c|c|c|c|c|}
\caption{Properties of new finite semifields of 64 elements}\label{superdatos}\\
\hline      \textbf{Plane}  & $\mathbf{S_3-class}$ & \textbf{$|$At$|$} & $\mathbf{(L_x,L_\infty,L_y)}$&\textbf{S/A sum}& $\mathbf{ZN}$\\\hline

\endfirsthead

\multicolumn{6}{c}%
{{ \tablename\ \thetable{} -- continued from previous page}} \\\hline \textbf{Plane}  & $\mathbf{S_3-class}$ & \textbf{$|$At$|$} & $\mathbf{(L_x,L_\infty,L_y)}$&\textbf{S/A sum}& $\mathbf{ZN}$\\ \hline
\endhead

\hline
\endfoot

\hline 
\endlastfoot

\textbf{	XIV	}		&$\begin{array}{c}\includegraphics[width=1cm,height=1cm]{hexagon.pdf}\\\end{array}$&	$  \begin{array}{c} 63 \\ \hbox{Solvable}\\\end{array}$ & $\begin{array}{c} 2[1]+1[63]\\ 3[1]+2[3]+2[7]+2[21]\\ 2[1]+1[63]\\\end{array}$	 & $\frac{63}{1}$	& $(2,2,2,2,4)$\\\hline	
\textbf{	XV	}		&$\begin{array}{c}\includegraphics[width=1cm,height=1cm]{hexagonLinTrB.pdf}\\\end{array}$&	$  \begin{array}{c} 21 \\ \hbox{Solvable}\\\end{array}$ & $\begin{array}{c} 2[1]+3[7]+2[21]\\ 2[1]+3[7]+2[21]\\ 2[1]+3[7]+2[21]\\\end{array}$	 & $\frac{186}{1} + \frac{9}{3}$	& $(2,2,2,2,2)$\\\hline	
\textbf{	XVI	}		&$\begin{array}{c}\includegraphics[width=1cm,height=1cm]{hexagonLinTrB.pdf}\\\end{array}$&	$  \begin{array}{c} 21 \\ \hbox{Solvable}\\\end{array}$ & $\begin{array}{c} 2[1]+3[7]+2[21]\\ 2[1]+3[7]+2[21]\\ 2[1]+3[7]+2[21]\\\end{array}$	 & $\frac{186}{1} + \frac{9}{3}$	& $(2,2,2,2,2)$\\\hline	
\textbf{	XVII	}		&$\begin{array}{c}\includegraphics[width=1cm,height=1cm]{hexagon.pdf}\\\end{array}$&	$  \begin{array}{c} 15 \\ \hbox{Cyclic}\\\end{array}$ & $\begin{array}{c} 2[1]+1[3]+4[15]\\ 2[1]+1[3]+4[15]\\ 5[1]+12[5]\\\end{array}$	 & $\frac{264}{1} + \frac{3}{5}$	& $(2,2,2,4,2)$\\\hline	
\textbf{	XVIII	}		&$\begin{array}{c}\includegraphics[width=1cm,height=1cm]{hexagon.pdf}\\\end{array}$&	$  \begin{array}{c} 15 \\ \hbox{Cyclic}\\\end{array}$ & $\begin{array}{c} 2[1]+1[3]+3[5]+3[15]\\ 2[1]+1[3]+4[15]\\ 2[1]+1[3]+4[15]\\\end{array}$	 & $\frac{264}{1} + \frac{3}{5}$	& $(2,2,2,2,2)$\\\hline	
\textbf{	XIX	}		&$\begin{array}{c}\includegraphics[width=1cm,height=1cm]{hexagonTriangTr.pdf}\\\end{array}$&	$  \begin{array}{c} 14 \\ \hbox{Cyclic}\\\end{array}$ & $\begin{array}{c} 2[1]+1[7]+4[14]\\ 2[1]+1[7]+4[14]\\ 2[1]+1[7]+4[14]\\\end{array}$	 & $\frac{280}{1} + \frac{7}{2}$	& $(2,2,2,2,2)$\\\hline	
\textbf{	XX	}		&$\begin{array}{c}\includegraphics[width=1cm,height=1cm]{hexagonLinTr.pdf}\\\end{array}$&	$  \begin{array}{c} 9 \\ \hbox{ $\ZZ_3^2 $}\\\end{array}$ & $\begin{array}{c} 2[1]+6[3]+5[9]\\ 2[1]+6[3]+5[9]\\ 2[1]+6[3]+5[9]\\\end{array}$	 & $\frac{440}{1} + \frac{3}{3}$	& $(2,2,2,2,2)$\\\hline	
\textbf{	XXI	}		&$\begin{array}{c}\includegraphics[width=1cm,height=1cm]{hexagon.pdf}\\\end{array}$&	$  \begin{array}{c} 9 \\ \hbox{ $\ZZ_3^2 $}\\\end{array}$ & $\begin{array}{c} 2[1]+6[3]+5[9]\\ 2[1]+6[3]+5[9]\\ 2[1]+6[3]+5[9]\\\end{array}$	 & $\frac{440}{1} + \frac{3}{3}$	& $(2,2,2,2,2)$\\\hline	
\textbf{	XXII	}		&$\begin{array}{c}\includegraphics[width=1cm,height=1cm]{hexagonLinTr.pdf}\\\end{array}$&	$  \begin{array}{c} 9 \\ \hbox{Cyclic}\\\end{array}$ & $\begin{array}{c} 2[1]+7[9]\\ 2[1]+7[9]\\ 5[1]+20[3]\\\end{array}$	 & $\frac{441}{1}$	& $(2,2,2,4,2)$\\\hline	
\textbf{	XXIII	}		&$\begin{array}{c}\includegraphics[width=1cm,height=1cm]{hexagonLinTrB.pdf}\\\end{array}$&	$  \begin{array}{c} 9 \\ \hbox{Cyclic}\\\end{array}$ & $\begin{array}{c} 5[1]+20[3]\\ 2[1]+7[9]\\ 2[1]+7[9]\\\end{array}$	 & $\frac{441}{1}$	& $(2,2,4,2,2)$\\\hline	
\textbf{	XXIV	}		&$\begin{array}{c}\includegraphics[width=1cm,height=1cm]{hexagonLinTr.pdf}\\\end{array}$&	$  \begin{array}{c} 9 \\ \hbox{ $\ZZ_3^2 $}\\\end{array}$ & $\begin{array}{c} 2[1]+21[3]\\ 2[1]+21[3]\\ 2[1]+6[3]+5[9]\\\end{array}$	 & $\frac{441}{1}$	& $(2,2,4,2,4)$\\\hline	
\textbf{	XXV	}		&$\begin{array}{c}\includegraphics[width=1cm,height=1cm]{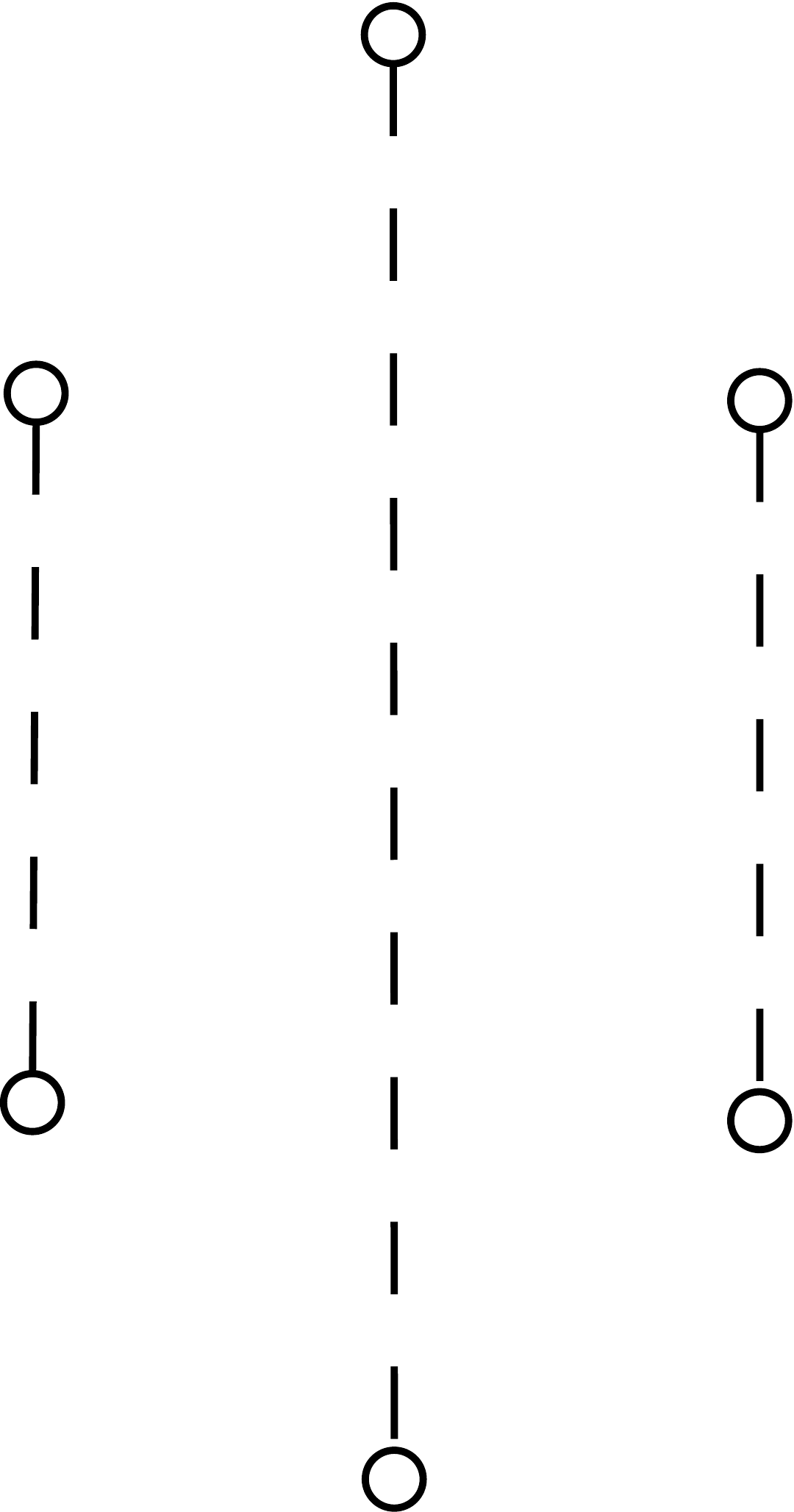}\\\end{array}$&	$  \begin{array}{c} 9 \\ \hbox{Cyclic}\\\end{array}$ & $\begin{array}{c} 2[1]+7[9]\\ 5[1]+20[3]\\ 2[1]+7[9]\\\end{array}$	 & $\frac{441}{1}$	& $(2,2,2,2,4)$\\\hline	
\textbf{	XXVI	}		&$\begin{array}{c}\includegraphics[width=1cm,height=1cm]{hexagonLinTr.pdf}\\\end{array}$&	$  \begin{array}{c} 9 \\ \hbox{ $\ZZ_3^2 $}\\\end{array}$ & $\begin{array}{c} 2[1]+21[3]\\ 2[1]+21[3]\\ 2[1]+6[3]+5[9]\\\end{array}$	 & $\frac{441}{1}$	& $(2,2,4,2,4)$\\\hline	
\textbf{	XXVII	}		&$\begin{array}{c}\includegraphics[width=1cm,height=1cm]{hexagonLinTr.pdf}\\\end{array}$&	$  \begin{array}{c} 9 \\ \hbox{ $\ZZ_3^2 $}\\\end{array}$ & $\begin{array}{c} 2[1]+21[3]\\ 2[1]+21[3]\\ 2[1]+6[3]+5[9]\\\end{array}$	 & $\frac{441}{1}$	& $(2,2,4,2,4)$\\\hline	
\textbf{	XXVIII	}		&$\begin{array}{c}\includegraphics[width=1cm,height=1cm]{hexagonLinTrC.pdf}\\\end{array}$&	$  \begin{array}{c} 9 \\ \hbox{ $\ZZ_3^2 $}\\\end{array}$ & $\begin{array}{c} 2[1]+6[3]+5[9]\\ 2[1]+6[3]+5[9]\\ 2[1]+6[3]+5[9]\\\end{array}$	 & $\frac{440}{1} + \frac{3}{3}$	& $(2,2,2,2,2)$\\\hline	
  
  \textbf{	XXIX	}		&$\begin{array}{c}\includegraphics[width=1cm,height=1cm]{hexagonTriangTr.pdf}\\\end{array}$&	$  \begin{array}{c} 9 \\ \hbox{ $\ZZ_3^2 $}\\\end{array}$ & $\begin{array}{c} 2[1]+6[3]+5[9]\\ 2[1]+6[3]+5[9]\\ 2[1]+6[3]+5[9]\\\end{array}$	 & $\frac{440}{1} + \frac{3}{3}$	& $(2,2,2,2,2)$\\\hline	

  \textbf{	XXX	}		&$\begin{array}{c}\includegraphics[width=1cm,height=1cm]{hexagonTriangTr.pdf}\\\end{array}$&	$  \begin{array}{c} 9 \\ \hbox{ $\ZZ_3^2 $}\\\end{array}$ & $\begin{array}{c} 2[1]+6[3]+5[9]\\ 2[1]+6[3]+5[9]\\ 2[1]+6[3]+5[9]\\\end{array}$	 & $\frac{440}{1} + \frac{3}{3}$	& $(2,2,2,2,2)$\\\hline	
\textbf{	XXXI	}		&$\begin{array}{c}\includegraphics[width=1cm,height=1cm]{hexagonLinTr.pdf}\\\end{array}$&	$  \begin{array}{c} 9 \\ \hbox{ $\ZZ_3^2 $}\\\end{array}$ & $\begin{array}{c} 2[1]+6[3]+5[9]\\ 2[1]+6[3]+5[9]\\ 2[1]+6[3]+5[9]\\\end{array}$	 & $\frac{440}{1} + \frac{3}{3}$	& $(2,2,2,2,2)$\\\hline	
\textbf{	XXXII	}		&$\begin{array}{c}\includegraphics[width=1cm,height=1cm]{hexagonCompTr.pdf}\\\end{array}$&	$  \begin{array}{c} 9 \\ \hbox{ $\ZZ_3^2 $}\\\end{array}$ & $\begin{array}{c} 2[1]+6[3]+5[9]\\ 2[1]+6[3]+5[9]\\ 2[1]+6[3]+5[9]\\\end{array}$	 & $\frac{440}{1} + \frac{3}{3}$	& $(2,2,2,2,2)$\\\hline	
\textbf{	XXXIII	}		&$\begin{array}{c}\includegraphics[width=1cm,height=1cm]{hexagon.pdf}\\\end{array}$&	$  \begin{array}{c} 7 \\ \hbox{Cyclic}\\\end{array}$ & $\begin{array}{c} 2[1]+9[7]\\ 9[1]+8[7]\\ 2[1]+9[7]\\\end{array}$	 & $\frac{567}{1}$	& $(2,2,2,2,2)$\\\hline	
\textbf{	XXXIV	}		&$\begin{array}{c}\includegraphics[width=1cm,height=1cm]{hexagon.pdf}\\\end{array}$&	$  \begin{array}{c} 6 \\ \hbox{Symmetric}\\\end{array}$ & $\begin{array}{c} 2[1]+7[3]+7[6]\\ 2[1]+7[3]+7[6]\\ 5[1]+6[2]+4[3]+6[6]\\\end{array}$	 & $\frac{637}{1} + \frac{49}{2}$	& $(2,2,2,2,2)$\\\hline	
\textbf{	XXXV	}		&$\begin{array}{c}\includegraphics[width=1cm,height=1cm]{hexagon.pdf}\\\end{array}$&	$  \begin{array}{c} 6 \\ \hbox{Cyclic}\\\end{array}$ & $\begin{array}{c} 3[1]+1[2]+2[3]+9[6]\\ 3[1]+1[2]+2[3]+9[6]\\ 3[1]+1[2]+2[3]+9[6]\\\end{array}$	 & $\frac{652}{1} + \frac{16}{2} + \frac{4}{3} + \frac{1}{6}$	& $(2,2,2,2,2)$\\\hline	
\textbf{	XXXVI	}		&$\begin{array}{c}\includegraphics[width=1cm,height=1cm]{hexagonLinTrB.pdf}\\\end{array}$&	$  \begin{array}{c} 6 \\ \hbox{Cyclic}\\\end{array}$ & $\begin{array}{c} 3[1]+1[2]+2[3]+9[6]\\ 3[1]+1[2]+2[3]+9[6]\\ 3[1]+1[2]+2[3]+9[6]\\\end{array}$	 & $\frac{652}{1} + \frac{16}{2} + \frac{4}{3} + \frac{1}{6}$	& $(2,2,2,2,2)$\\\hline	
\textbf{	XXXVII	}		&$\begin{array}{c}\includegraphics[width=1cm,height=1cm]{hexagon.pdf}\\\end{array}$&	$  \begin{array}{c} 6 \\ \hbox{Symmetric}\\\end{array}$ & $\begin{array}{c} 2[1]+7[3]+7[6]\\ 2[1]+7[3]+7[6]\\ 5[1]+6[2]+4[3]+6[6]\\\end{array}$	 & $\frac{637}{1} + \frac{49}{2}$	& $(2,2,2,2,2)$\\\hline	
\textbf{	XXXVIII	}		&$\begin{array}{c}\includegraphics[width=1cm,height=1cm]{hexagon.pdf}\\\end{array}$&	$  \begin{array}{c} 6 \\ \hbox{Symmetric}\\\end{array}$ & $\begin{array}{c} 2[1]+7[3]+7[6]\\ 2[1]+7[3]+7[6]\\ 5[1]+6[2]+4[3]+6[6]\\\end{array}$	 & $\frac{637}{1} + \frac{49}{2}$	& $(2,2,2,2,2)$\\\hline	
\textbf{	XXXIX	}		&$\begin{array}{c}\includegraphics[width=1cm,height=1cm]{hexagonLinTrB.pdf}\\\end{array}$&	$  \begin{array}{c} 5 \\ \hbox{Cyclic}\\\end{array}$ & $\begin{array}{c} 5[1]+12[5]\\ 5[1]+12[5]\\ 5[1]+12[5]\\\end{array}$	 & $\frac{792}{1} + \frac{9}{5}$	& $(2,2,2,2,2)$\\\hline	
\textbf{	XL	}		&$\begin{array}{c}\includegraphics[width=1cm,height=1cm]{hexagon.pdf}\\\end{array}$&	$  \begin{array}{c} 3 \\ \hbox{Cyclic}\\\end{array}$ & $\begin{array}{c} 2[1]+21[3]\\ 2[1]+21[3]\\ 65[1]\\\end{array}$	 & $\frac{1323}{1}$	& $(2,2,2,4,2)$\\\hline	
\textbf{	XLI	}		&$\begin{array}{c}\includegraphics[width=1cm,height=1cm]{hexagon.pdf}\\\end{array}$&	$  \begin{array}{c} 3 \\ \hbox{Cyclic}\\\end{array}$ & $\begin{array}{c} 65[1]\\ 2[1]+21[3]\\ 2[1]+21[3]\\\end{array}$	 & $\frac{1323}{1}$	& $(2,2,4,2,2)$\\\hline	
\textbf{	XLII	}		&$\begin{array}{c}\includegraphics[width=1cm,height=1cm]{hexagon.pdf}\\\end{array}$&	$  \begin{array}{c} 3 \\ \hbox{Cyclic}\\\end{array}$ & $\begin{array}{c} 2[1]+21[3]\\ 65[1]\\ 2[1]+21[3]\\\end{array}$	 & $\frac{1323}{1}$	& $(2,2,2,2,4)$\\\hline	
\textbf{	XLIII	}		&$\begin{array}{c}\includegraphics[width=1cm,height=1cm]{hexagon.pdf}\\\end{array}$&	$  \begin{array}{c} 3 \\ \hbox{Cyclic}\\\end{array}$ & $\begin{array}{c} 65[1]\\ 2[1]+21[3]\\ 2[1]+21[3]\\\end{array}$	 & $\frac{1323}{1}$	& $(2,2,4,2,2)$\\\hline	  \textbf{	XLIV	}		&$\begin{array}{c}\includegraphics[width=1cm,height=1cm]{hexagon.pdf}\\\end{array}$&	$  \begin{array}{c} 3 \\ \hbox{Cyclic}\\\end{array}$ & $\begin{array}{c} 2[1]+21[3]\\ 2[1]+21[3]\\ 65[1]\\\end{array}$	 & $\frac{1323}{1}$	& $(2,2,2,4,2)$\\\hline	
\textbf{	XLV	}		&$\begin{array}{c}\includegraphics[width=1cm,height=1cm]{hexagon.pdf}\\\end{array}$&	$  \begin{array}{c} 3 \\ \hbox{Cyclic}\\\end{array}$ & $\begin{array}{c} 65[1]\\ 2[1]+21[3]\\ 2[1]+21[3]\\\end{array}$	 & $\frac{1323}{1}$	& $(2,2,4,2,2)$\\\hline	
  \textbf{	XLVI	}		&$\begin{array}{c}\includegraphics[width=1cm,height=1cm]{hexagon.pdf}\\\end{array}$&	$  \begin{array}{c} 3 \\ \hbox{Cyclic}\\\end{array}$ & $\begin{array}{c} 2[1]+21[3]\\ 17[1]+16[3]\\ 2[1]+21[3]\\\end{array}$	 & $\frac{1323}{1}$	& $(2,2,2,2,2)$\\\hline	
\textbf{	XLVII	}		&$\begin{array}{c}\includegraphics[width=1cm,height=1cm]{hexagon.pdf}\\\end{array}$&	$  \begin{array}{c} 3 \\ \hbox{Cyclic}\\\end{array}$ & $\begin{array}{c} 2[1]+21[3]\\ 2[1]+21[3]\\ 17[1]+16[3]\\\end{array}$	 & $\frac{1323}{1}$	& $(2,2,2,2,2)$\\\hline	
\textbf{	XLVIII	}		&$\begin{array}{c}\includegraphics[width=1cm,height=1cm]{hexagonLinTr.pdf}\\\end{array}$&	$  \begin{array}{c} 3 \\ \hbox{Cyclic}\\\end{array}$ & $\begin{array}{c} 2[1]+21[3]\\ 2[1]+21[3]\\ 17[1]+16[3]\\\end{array}$	 & $\frac{1323}{1}$	& $(2,2,2,2,2)$\\\hline	
\textbf{	XLIX	}		&$\begin{array}{c}\includegraphics[width=1cm,height=1cm]{hexagonLinTrC.pdf}\\\end{array}$&	$  \begin{array}{c} 3 \\ \hbox{Cyclic}\\\end{array}$ & $\begin{array}{c} 2[1]+21[3]\\ 17[1]+16[3]\\ 2[1]+21[3]\\\end{array}$	 & $\frac{1323}{1}$	& $(2,2,2,2,2)$\\\hline	
\textbf{	L	}		&$\begin{array}{c}\includegraphics[width=1cm,height=1cm]{hexagonLinTr.pdf}\\\end{array}$&	$  \begin{array}{c} 3 \\ \hbox{Cyclic}\\\end{array}$ & $\begin{array}{c} 2[1]+21[3]\\ 2[1]+21[3]\\ 17[1]+16[3]\\\end{array}$	 & $\frac{1323}{1}$	& $(2,2,2,2,2)$\\\hline	
\textbf{	LI	}		&$\begin{array}{c}\includegraphics[width=1cm,height=1cm]{hexagonLinTrB.pdf}\\\end{array}$&	$  \begin{array}{c} 3 \\ \hbox{Cyclic}\\\end{array}$ & $\begin{array}{c} 17[1]+16[3]\\ 2[1]+21[3]\\ 2[1]+21[3]\\\end{array}$	 & $\frac{1323}{1}$	& $(2,2,2,2,2)$\\\hline	
\textbf{	LII	}		&$\begin{array}{c}\includegraphics[width=1cm,height=1cm]{hexagonLinTrB.pdf}\\\end{array}$&	$  \begin{array}{c} 3 \\ \hbox{Cyclic}\\\end{array}$ & $\begin{array}{c} 17[1]+16[3]\\ 2[1]+21[3]\\ 2[1]+21[3]\\\end{array}$	 & $\frac{1323}{1}$	& $(2,2,2,2,2)$\\\hline	
\textbf{	LIII	}		&$\begin{array}{c}\includegraphics[width=1cm,height=1cm]{hexagon.pdf}\\\end{array}$&	$  \begin{array}{c} 3 \\ \hbox{Cyclic}\\\end{array}$ & $\begin{array}{c} 65[1]\\ 2[1]+21[3]\\ 2[1]+21[3]\\\end{array}$	 & $\frac{1323}{1}$	& $(2,2,4,2,2)$\\\hline	
\textbf{	LIV	}		&$\begin{array}{c}\includegraphics[width=1cm,height=1cm]{hexagon.pdf}\\\end{array}$&	$  \begin{array}{c} 3 \\ \hbox{Cyclic}\\\end{array}$ & $\begin{array}{c} 65[1]\\ 2[1]+21[3]\\ 2[1]+21[3]\\\end{array}$	 & $\frac{1323}{1}$	& $(2,2,4,2,2)$\\\hline	
\textbf{	LV	}		&$\begin{array}{c}\includegraphics[width=1cm,height=1cm]{hexagon.pdf}\\\end{array}$&	$  \begin{array}{c} 3 \\ \hbox{Cyclic}\\\end{array}$ & $\begin{array}{c} 2[1]+21[3]\\ 65[1]\\ 2[1]+21[3]\\\end{array}$	 & $\frac{1323}{1}$	& $(2,2,2,2,4)$\\\hline	
\textbf{	LVI	}		&$\begin{array}{c}\includegraphics[width=1cm,height=1cm]{hexagonLinTrB.pdf}\\\end{array}$&	$  \begin{array}{c} 3 \\ \hbox{Cyclic}\\\end{array}$ & $\begin{array}{c} 17[1]+16[3]\\ 2[1]+21[3]\\ 2[1]+21[3]\\\end{array}$	 & $\frac{1323}{1}$	& $(2,2,2,2,2)$\\\hline	
\textbf{	LVII	}		&$\begin{array}{c}\includegraphics[width=1cm,height=1cm]{hexagon.pdf}\\\end{array}$&	$  \begin{array}{c} 3 \\ \hbox{Cyclic}\\\end{array}$ & $\begin{array}{c} 17[1]+16[3]\\ 2[1]+21[3]\\ 2[1]+21[3]\\\end{array}$	 & $\frac{1323}{1}$	& $(2,2,2,2,2)$\\\hline	
\textbf{	LVIII	}		&$\begin{array}{c}\includegraphics[width=1cm,height=1cm]{hexagon.pdf}\\\end{array}$&	$  \begin{array}{c} 3 \\ \hbox{Cyclic}\\\end{array}$ & $\begin{array}{c} 65[1]\\ 2[1]+21[3]\\ 2[1]+21[3]\\\end{array}$	 & $\frac{1323}{1}$	& $(2,2,4,2,2)$\\\hline	
  	\textbf{	LIX	}		&$\begin{array}{c}\includegraphics[width=1cm,height=1cm]{hexagon.pdf}\\\end{array}$&	$  \begin{array}{c} 3 \\ \hbox{Cyclic}\\\end{array}$ & $\begin{array}{c} 17[1]+16[3]\\ 2[1]+21[3]\\ 2[1]+21[3]\\\end{array}$	 & $\frac{1323}{1}$	& $(2,2,2,2,2)$\\\hline	
\textbf{	LX	}		&$\begin{array}{c}\includegraphics[width=1cm,height=1cm]{hexagonLinTrC.pdf}\\\end{array}$&	$  \begin{array}{c} 3 \\ \hbox{Cyclic}\\\end{array}$ & $\begin{array}{c} 2[1]+21[3]\\ 17[1]+16[3]\\ 2[1]+21[3]\\\end{array}$	 & $\frac{1323}{1}$	& $(2,2,2,2,2)$\\\hline	
\textbf{	LXI	}		&$\begin{array}{c}\includegraphics[width=1cm,height=1cm]{hexagon.pdf}\\\end{array}$&	$  \begin{array}{c} 3 \\ \hbox{Cyclic}\\\end{array}$ & $\begin{array}{c} 17[1]+16[3]\\ 2[1]+21[3]\\ 2[1]+21[3]\\\end{array}$	 & $\frac{1323}{1}$	& $(2,2,2,2,2)$\\\hline
\textbf{	LXII	}		&$\begin{array}{c}\includegraphics[width=1cm,height=1cm]{hexagonLinTr.pdf}\\\end{array}$&	$  \begin{array}{c} 3 \\ \hbox{Cyclic}\\\end{array}$ & $\begin{array}{c} 2[1]+21[3]\\ 2[1]+21[3]\\ 17[1]+16[3]\\\end{array}$	 & $\frac{1323}{1}$	& $(2,2,2,2,2)$\\\hline	
\textbf{	LXIII	}		&$\begin{array}{c}\includegraphics[width=1cm,height=1cm]{hexagon.pdf}\\\end{array}$&	$  \begin{array}{c} 3 \\ \hbox{Cyclic}\\\end{array}$ & $\begin{array}{c} 65[1]\\ 2[1]+21[3]\\ 2[1]+21[3]\\\end{array}$	 & $\frac{1323}{1}$	& $(2,2,4,2,2)$\\\hline	
\textbf{	LXIV	}		&$\begin{array}{c}\includegraphics[width=1cm,height=1cm]{hexagonLinTr.pdf}\\\end{array}$&	$  \begin{array}{c} 3 \\ \hbox{Cyclic}\\\end{array}$ & $\begin{array}{c} 2[1]+21[3]\\ 2[1]+21[3]\\ 17[1]+16[3]\\\end{array}$	 & $\frac{1323}{1}$	& $(2,2,2,2,2)$\\\hline	
\textbf{	LXV	}		&$\begin{array}{c}\includegraphics[width=1cm,height=1cm]{hexagon.pdf}\\\end{array}$&	$  \begin{array}{c} 3 \\ \hbox{Cyclic}\\\end{array}$ & $\begin{array}{c} 17[1]+16[3]\\ 2[1]+21[3]\\ 2[1]+21[3]\\\end{array}$	 & $\frac{1323}{1}$	& $(2,2,2,2,2)$\\\hline	
\textbf{	LXVI	}		&$\begin{array}{c}\includegraphics[width=1cm,height=1cm]{hexagonLinTrC.pdf}\\\end{array}$&	$  \begin{array}{c} 3 \\ \hbox{Cyclic}\\\end{array}$ & $\begin{array}{c} 2[1]+21[3]\\ 17[1]+16[3]\\ 2[1]+21[3]\\\end{array}$	 & $\frac{1323}{1}$	& $(2,2,2,2,2)$\\\hline	
\textbf{	LXVII	}		&$\begin{array}{c}\includegraphics[width=1cm,height=1cm]{hexagon.pdf}\\\end{array}$&	$  \begin{array}{c} 3 \\ \hbox{Cyclic}\\\end{array}$ & $\begin{array}{c} 65[1]\\ 2[1]+21[3]\\ 2[1]+21[3]\\\end{array}$	 & $\frac{1323}{1}$	& $(2,2,4,2,2)$\\\hline	
\textbf{	LXVIII	}		&$\begin{array}{c}\includegraphics[width=1cm,height=1cm]{hexagon.pdf}\\\end{array}$&	$  \begin{array}{c} 3 \\ \hbox{Cyclic}\\\end{array}$ & $\begin{array}{c} 65[1]\\ 2[1]+21[3]\\ 2[1]+21[3]\\\end{array}$	 & $\frac{1323}{1}$	& $(2,2,4,2,2)$\\\hline	
\textbf{	LXIX	}		&$\begin{array}{c}\includegraphics[width=1cm,height=1cm]{hexagon.pdf}\\\end{array}$&	$  \begin{array}{c} 3 \\ \hbox{Cyclic}\\\end{array}$ & $\begin{array}{c} 65[1]\\ 2[1]+21[3]\\ 2[1]+21[3]\\\end{array}$	 & $\frac{1323}{1}$	& $(2,2,4,2,2)$\\\hline	
\textbf{	LXX	}		&$\begin{array}{c}\includegraphics[width=1cm,height=1cm]{hexagon.pdf}\\\end{array}$&	$  \begin{array}{c} 3 \\ \hbox{Cyclic}\\\end{array}$ & $\begin{array}{c} 2[1]+21[3]\\ 2[1]+21[3]\\ 65[1]\\\end{array}$	 & $\frac{1323}{1}$	& $(2,2,2,4,2)$\\\hline	
\textbf{	LXXI	}		&$\begin{array}{c}\includegraphics[width=1cm,height=1cm]{hexagon.pdf}\\\end{array}$&	$  \begin{array}{c} 2 \\ \hbox{Cyclic}\\\end{array}$ & $\begin{array}{c} 9[1]+28[2]\\ 9[1]+28[2]\\ 9[1]+28[2]\\\end{array}$	 & $\frac{1960}{1} + \frac{49}{2}$	& $(2,2,2,2,2)$\\\hline	
\textbf{	LXXII	}		&$\begin{array}{c}\includegraphics[width=1cm,height=1cm]{hexagon.pdf}\\\end{array}$&	$  \begin{array}{c} 2 \\ \hbox{Cyclic}\\\end{array}$ & $\begin{array}{c} 9[1]+28[2]\\ 9[1]+28[2]\\ 9[1]+28[2]\\\end{array}$	 & $\frac{1960}{1} + \frac{49}{2}$	& $(2,2,2,2,2)$\\\hline	
\textbf{	LXXIII	}		&$\begin{array}{c}\includegraphics[width=1cm,height=1cm]{hexagon.pdf}\\\end{array}$&	$  \begin{array}{c} 2 \\ \hbox{Cyclic}\\\end{array}$ & $\begin{array}{c} 9[1]+28[2]\\ 9[1]+28[2]\\ 9[1]+28[2]\\\end{array}$	 & $\frac{1960}{1} + \frac{49}{2}$	& $(2,2,2,2,2)$\\\hline
\textbf{	LXXIV	}		&$\begin{array}{c}\includegraphics[width=1cm,height=1cm]{hexagon.pdf}\\\end{array}$&	$  \begin{array}{c} 2 \\ \hbox{Cyclic}\\\end{array}$ & $\begin{array}{c} 9[1]+28[2]\\ 9[1]+28[2]\\ 9[1]+28[2]\\\end{array}$	 & $\frac{1960}{1} + \frac{49}{2}$	& $(2,2,2,2,2)$\\\hline	
\textbf{	LXXV	}		&$\begin{array}{c}\includegraphics[width=1cm,height=1cm]{hexagonTriangTr.pdf}\\\end{array}$&	$  \begin{array}{c} 1 \\ \hbox{Cyclic}\\\end{array}$ & $\begin{array}{c} 65[1]\\ 65[1]\\ 65[1]\\\end{array}$	 & $\frac{3969}{1}$	& $(2,2,2,2,2)$\\\hline	
\textbf{	LXXVI	}		&$\begin{array}{c}\includegraphics[width=1cm,height=1cm]{hexagonLinTrC.pdf}\\\end{array}$&	$  \begin{array}{c} 1 \\ \hbox{Cyclic}\\\end{array}$ & $\begin{array}{c} 65[1]\\ 65[1]\\ 65[1]\\\end{array}$	 & $\frac{3969}{1}$	& $(2,2,2,2,2)$\\\hline	
\textbf{	LXXVII	}		&$\begin{array}{c}\includegraphics[width=1cm,height=1cm]{hexagonTriangTr.pdf}\\\end{array}$&	$  \begin{array}{c} 1 \\ \hbox{Cyclic}\\\end{array}$ & $\begin{array}{c} 65[1]\\ 65[1]\\ 65[1]\\\end{array}$	 & $\frac{3969}{1}$	& $(2,2,2,2,2)$\\\hline
  \textbf{	LXXVIII	}		&$\begin{array}{c}\includegraphics[width=1cm,height=1cm]{hexagonTriangTr.pdf}\\\end{array}$&	$  \begin{array}{c} 1 \\ \hbox{Cyclic}\\\end{array}$ & $\begin{array}{c} 65[1]\\ 65[1]\\ 65[1]\\\end{array}$	 & $\frac{3969}{1}$	& $(2,2,2,2,2)$\\\hline	
\textbf{	LXXIX	}		&$\begin{array}{c}\includegraphics[width=1cm,height=1cm]{hexagonTriangTr.pdf}\\\end{array}$&	$  \begin{array}{c} 1 \\ \hbox{Cyclic}\\\end{array}$ & $\begin{array}{c} 65[1]\\ 65[1]\\ 65[1]\\\end{array}$	 & $\frac{3969}{1}$	& $(2,2,2,2,2)$\\\hline	
\textbf{	LXXX	}		&$\begin{array}{c}\includegraphics[width=1cm,height=1cm]{hexagon.pdf}\\\end{array}$&	$  \begin{array}{c} 1 \\ \hbox{Cyclic}\\\end{array}$ & $\begin{array}{c} 65[1]\\ 65[1]\\ 65[1]\\\end{array}$	 & $\frac{3969}{1}$	& $(2,2,2,2,2)$\\\hline

\end{longtable}
\end{center}

The autotopism group of the semifield XIV is $<x,y,z\ |\ x^3=y,y^3=z^7=1,z^x=z^2>$. The autotopism group of the semifields XV and XVI is $<x,y\ |\ x^3=y^7=1,y^x=y^2>$.

\section{Concluding remarks}

In this paper we present a classification of finite semifields of order 64, with the help of computational tools. The resulting finite semifield \emph{zoo} is very rich, and it completes the classification of finite semifields up to order 125.
The algorithm used is a combination of known methods and some crucial observations. It proved to be quite suitable to the problem considered. We hope this classification will be helpful in the discovery of new semifield constructions, a problem which has not been addressed in this paper.

\section*{Appendix}

Using the Oyama representation of matrices \cite{Oyama} we managed to obtain different algebraic descriptions for the semifield representatives. As a complement we include these presentations here, in the hope that they can be helpful in the discovery of new constructions. There are three different types of algebraic representations:

\begin{enumerate}
	\item $D=GF(8)^2$, with $j^3+j+1=0$, and multiplication rule $(x_1,x_2)*(a_1,a_2)$.
	\item $D=GF(4)^3$, with $j^2+j+1=0$,  and multiplication rule $(x_1,x_2,x_3)*(a_1,a_2,a_3)$.
	\item $D=GF(64)$, with $j^6+j+1=0$, and multiplication rule $x*a$.
\end{enumerate}

Next we list the multiplication rules of the semifield representatives.

\begin{itemize}

\item{\bf	XIV:	(1)}	$		    (x_1 a_1  + x_2 (j^{4} a_1+j^{4} a_1^{4} + a_2^{4})  + x_2^{2} (j^{3} a_1+j^{3} a_1^{4} + j^{3} a_2^{4})  + x_2^{4} (j^{6} a_1+j^{6} a_1^{4} + j^{6} a_2^{4}) \ ,\  x_1 a_2  + x_2 (a_1^{4} + j^{4} a_2)  + x_2^{2} (j^{3} a_2+a_2^{4})  + x_2^{4} (j^{6} a_2+j^{3} a_2^{4}) )		$
\item{\bf	XV:	(1)}	$		    (x_1 a_1  + x_2 (j a_2+j^{6} a_2^{4})  + x_2^{2} j^{2} a_2^{2} + x_2^{4} (j^{2} a_1+j^{2} a_1^{2} + j^{6} a_2+j^{3} a_2^{2}) \ ,\  x_1 a_2  + x_2 (a_1 + j^{2} a_2^{4})  + x_2^{2} j^{5} a_2^{2} + x_2^{4} j^{2} a_2)		$
\item{\bf	XVI:	(1)}	$		    (x_1 a_1  + x_2 (j^{2} a_2+a_2^{2})  + x_2^{2} (j^{6} a_1+j^{6} a_1^{4} + j^{3} a_2^{4})  + x_2^{4} (j^{5} a_1+j^{5} a_1^{2} + a_2^{2}+j^{3} a_2^{4}) \ ,\  x_1 a_2  + x_2 (a_1 + j^{6} a_2^{2})  + x_2^{2} j^{6} a_2 + x_2^{4} j^{5} a_2)		$
\item{\bf	XVII:	(2)}	$		    (x_1 a_1  + x_2 (j^{2} a_2 + j^{2} a_3)  + x_2^{2} j^{2} a_2^{2} + x_3 (j a_2 + j^{2} a_3)  + x_3^{2} (j^{2} a_2^{2} + j^{2} a_3^{2}) \ ,\  x_1 a_2  + x_2 a_1 + x_2^{2} j^{2} a_3^{2} + x_3^{2} (j^{2} a_2^{2} + j^{2} a_3^{2}) \ ,\  x_1 a_3  + x_2 j^{2} a_3 + x_2^{2} (a_2^{2} + j^{2} a_3^{2})  + x_3 (a_1 + j^{2} a_2)  + x_3^{2} (j^{2} a_2^{2} + a_3^{2}) )		$
\item{\bf	XVIII: (2)	}	$		    (x_1 a_1  + x_2 (j^{2} a_1+j^{2} a_1^{2} + a_2^{2} + j^{2} a_3+j a_3^{2})  + x_2^{2} (a_2^{2} + j^{2} a_3)  + x_3 (j a_1+j a_1^{2} + j^{2} a_2+j^{2} a_2^{2} + a_3^{2})  + x_3^{2} (a_2 + a_3+j^{2} a_3^{2}) \ ,\  x_1 a_2  + x_2 (a_1^{2} + j^{2} a_2+j a_2^{2} + j^{2} a_3^{2})  + x_2^{2} (j^{2} a_2+j a_2^{2} + j^{2} a_3+j a_3^{2})  + x_3 j a_2 + x_3^{2} (j a_2+j a_2^{2} + j^{2} a_3) \ ,\  x_1 a_3  + x_2 j^{2} a_3 + x_2^{2} (a_2 + j^{2} a_3+j^{2} a_3^{2})  + x_3 (a_1^{2} + j a_2^{2} + j a_3+j^{2} a_3^{2})  + x_3^{2} (j^{2} a_2^{2} + j a_3+j a_3^{2}) )		$
\item{\bf	XIX:	(3)}	$			    (x_1 (j^{59} a_1 + j^{7} a_1^{2} + j^{30} a_1^{4} + j^{8} a_1^{8} + j^{11} a_1^{16} + j^{17} a_1^{32}) + x_1^{2} (j^{5} a_1 + j^{39} a_1^{2} + j^{6} a_1^{4} + j^{13} a_1^{8} + j^{57} a_1^{16} + j^{62} a_1^{32}) + x_1^{4} (j^{4} a_1 + j^{7} a_1^{2} + j^{9} a_1^{4} + j^{53} a_1^{8} + j^{49} a_1^{16} + j^{57} a_1^{32}) + x_1^{8} (j^{29} a_1 + j^{37} a_1^{2} + j^{22} a_1^{4} + j^{23} a_1^{8} + j^{51} a_1^{16} + j^{11} a_1^{32}) + x_1^{16} (j^{16} a_1 + j^{12} a_1^{2} + j^{43} a_1^{4} + j^{16} a_1^{8} + j^{51} a_1^{16} + j^{45} a_1^{32}) + x_1^{32} (j^{47} a_1 + j^{56} a_1^{2} + j^{53} a_1^{4} + j^{36} a_1^{8} + j^{26} a_1^{16} + j^{6} a_1^{32}))		$
\item{\bf	XX:	(2)}	$		    (x_1 a_1  + x_2 (a_2^{2} + a_3^{2})  + x_2^{2} (a_2^{2} + j a_3^{2})  + x_3 (a_1+a_1^{2} + a_2^{2} + a_3^{2})  + x_3^{2} (j a_2^{2} + a_3^{2}) \ ,\  x_1 a_2  + x_2 (a_1^{2} + a_3^{2})  + x_3 a_2 + x_3^{2} j a_3^{2}\ ,\  x_1 a_3  + x_2^{2} (a_2^{2} + j a_3^{2})  + x_3 (a_1^{2} + a_3+a_3^{2})  + x_3^{2} (j a_2^{2} + a_3^{2}) )		$
\item{\bf	XXI:	(2)}	$		    (x_1 a_1  + x_2 (a_1+a_1^{2} + j a_3^{2})  + x_3 j a_3^{2} + x_3^{2} (j^{2} a_2^{2} + j^{2} a_3^{2}) \ ,\  x_1 a_2  + x_2 (a_1^{2} + a_2+a_2^{2})  + x_3^{2} (j^{2} a_2^{2} + j^{2} a_3^{2}) \ ,\  x_1 a_3  + x_2 a_3 + x_2^{2} a_2^{2} + x_3 (a_1^{2} + a_2^{2})  + x_3^{2} (j a_2^{2} + a_3^{2}) )		$
\item{\bf	XXII:	(2)}	$		    (x_1 a_1  + x_2 (j^{2} a_1+j^{2} a_1^{2} + a_2^{2} + j^{2} a_3)  + x_2^{2} (j^{2} a_2+j a_2^{2} + j a_3+a_3^{2})  + x_3 j^{2} a_3 + x_3^{2} (j a_2^{2} + a_3^{2}) \ ,\  x_1 a_2  + x_2 (a_1^{2} + j^{2} a_3^{2})  + x_2^{2} (j a_2+j a_2^{2} + j a_3+j^{2} a_3^{2})  + x_3 j^{2} a_2 + x_3^{2} (j a_2^{2} + j^{2} a_3^{2}) \ ,\  x_1 a_3  + x_2 (a_1+a_1^{2} + j a_2+j a_2^{2} + a_3)  + x_2^{2} (j a_2+j^{2} a_2^{2} + j a_3+j a_3^{2})  + x_3 (a_1 + j a_2 + j a_3)  + x_3^{2} (j^{2} a_2^{2} + j a_3^{2}) )		$
\item{\bf	XXIII: (2)	}	$		    (x_1 a_1  + x_2 (j a_1+j a_1^{2} + j a_2+j a_2^{2} + j a_3^{2})  + x_3 (a_2 + j a_3+j a_3^{2}) \ ,\  x_1 a_2  + x_2 (j a_1+j^{2} a_1^{2} + j^{2} a_2+j^{2} a_2^{2} + j a_3+j a_3^{2})  + x_3 (j a_1+j a_1^{2} + a_2+j^{2} a_2^{2} + a_3) \ ,\  x_1 a_3  + x_2 (j^{2} a_1+j^{2} a_1^{2} + j^{2} a_2+j a_2^{2} + j^{2} a_3+j a_3^{2})  + x_3 (a_1^{2} + a_2 + j^{2} a_3^{2}) )		$
\item{\bf	XXIV: (2)	}	$		    (x_1 a_1  + x_2 (j a_1+j a_1^{2} + j^{2} a_2+j^{2} a_2^{2} + j a_3+j^{2} a_3^{2})  + x_3 (a_1+a_1^{2} + j^{2} a_2+a_2^{2} + j^{2} a_3^{2}) \ ,\  x_1 a_2  + x_2 (j^{2} a_1+j a_1^{2} + a_2+a_2^{2} + a_3+j^{2} a_3^{2})  + x_3 (a_1+a_1^{2} + a_2+a_2^{2} + j^{2} a_3+a_3^{2}) \ ,\  x_1 a_3  + x_2 (j a_1+j a_1^{2} + j^{2} a_2+j a_2^{2} + a_3+j^{2} a_3^{2})  + x_3 (a_1 + j a_2 + j a_3) )		$
\item{\bf	XXV:	 (1)}	$		    (x_1 a_1  + x_2 (j^{5} a_1+j^{3} a_1^{2}+j^{2} a_1^{4} + a_2+j a_2^{2}+j^{5} a_2^{4})  + x_2^{2} (j^{2} a_1+a_1^{2}+j^{6} a_1^{4} + j^{5} a_2+j a_2^{2})  + x_2^{4} (j^{4} a_1+a_1^{2}+j^{5} a_1^{4} + j^{3} a_2+j^{5} a_2^{4}) \ ,\  x_1 a_2  + x_2 (j a_1+j^{6} a_1^{2}+j^{4} a_1^{4} + a_2+j a_2^{2}+j a_2^{4})  + x_2^{2} (j^{6} a_1+j^{6} a_1^{2} + j^{2} a_2^{2}+j a_2^{4})  + x_2^{4} (j^{4} a_1+j^{4} a_1^{4} + j^{5} a_2+j^{3} a_2^{2}+j^{6} a_2^{4}) )		$
\item{\bf	XXVI: (2)	}	$		    (x_1 a_1  + x_2 (a_1+a_1^{2} + a_2+a_2^{2} + j a_3)  + x_3 (a_1+a_1^{2} + j a_2 + a_3+j^{2} a_3^{2}) \ ,\  x_1 a_2  + x_2 (a_1 + j^{2} a_3)  + x_3 (a_1+a_1^{2} + j^{2} a_2+j a_2^{2} + j^{2} a_3+a_3^{2}) \ ,\  x_1 a_3  + x_2 (a_2^{2} + j a_3+j a_3^{2})  + x_3 (a_1 + j a_2+j^{2} a_2^{2} + j^{2} a_3^{2}) )		$
\item{\bf	XXVII: (2)	}	$		    (x_1 a_1  + x_2 (j a_1+j a_1^{2} + j a_2+j a_2^{2} + a_3+a_3^{2})  + x_3 (j a_1+j a_1^{2} + j^{2} a_2+a_2^{2} + j^{2} a_3) \ ,\  x_1 a_2  + x_2 (a_1 + a_2+a_2^{2} + a_3)  + x_3 (j a_1+j a_1^{2} + j^{2} a_2+j^{2} a_2^{2} + j a_3^{2}) \ ,\  x_1 a_3  + x_2 (j^{2} a_1+j^{2} a_1^{2} + a_2 + j a_3^{2})  + x_3 (a_1 + j a_2 + j a_3) )		$
\item{\bf	XXVIII: (2)	}	$		    (x_1 a_1  + x_2 j a_2^{2} + x_2^{2} (j a_2^{2} + j a_3^{2})  + x_3 (j a_1+j a_1^{2} + j a_2^{2} + j a_3^{2}) \ ,\  x_1 a_2  + x_2 (a_1^{2} + j^{2} a_3^{2})  + x_2^{2} j a_3^{2} + x_3 j a_2 + x_3^{2} (j^{2} a_2^{2} + j^{2} a_3^{2}) \ ,\  x_1 a_3  + x_2^{2} (a_2^{2} + a_3^{2})  + x_3 (a_1^{2} + j a_3+j^{2} a_3^{2}) )		$
\item{\bf	XXIX: (2)}	$		    (x_1 a_1  + x_2 (a_1+a_1^{2} + a_2^{2} + j^{2} a_3^{2})  + x_2^{2} (a_2^{2} + a_3^{2})  + x_3 (a_1+a_1^{2} + j^{2} a_2^{2} + a_3^{2})  + x_3^{2} a_2^{2}\ ,\  x_1 a_2  + x_2 (a_1^{2} + a_2+a_2^{2} + a_3^{2})  + x_2^{2} j^{2} a_3^{2} + x_3 a_2 + x_3^{2} (a_2^{2} + a_3^{2}) \ ,\  x_1 a_3  + x_2 a_3 + x_2^{2} (a_2^{2} + j a_3^{2})  + x_3 (a_1^{2} + a_2^{2} + a_3+a_3^{2})  + x_3^{2} a_3^{2})		$
\item{\bf	XXX:	 (2)}	$		    (x_1 a_1  + x_2 (a_1+a_1^{2})  + x_2^{2} j a_3^{2} + x_3 (a_2^{2} + j^{2} a_3^{2})  + x_3^{2} a_3^{2}\ ,\  x_1 a_2  + x_2 (a_1^{2} + a_2+a_2^{2})  + x_2^{2} j a_3^{2} + x_3^{2} a_3^{2}\ ,\  x_1 a_3  + x_2 a_3 + x_2^{2} (a_2^{2} + j^{2} a_3^{2})  + x_3 (a_1^{2} + a_2^{2})  + x_3^{2} (j^{2} a_2^{2} + j^{2} a_3^{2}) )		$
\item{\bf	XXXI: (2)	}	$		    (x_1 a_1  + x_2 (a_1+a_1^{2} + j^{2} a_2^{2} + j^{2} a_3^{2})  + x_2^{2} (j^{2} a_2^{2} + a_3^{2})  + x_3 (j^{2} a_1+j^{2} a_1^{2} + j a_2^{2})  + x_3^{2} a_3^{2}\ ,\  x_1 a_2  + x_2 (a_1^{2} + a_2+a_2^{2} + j a_3^{2})  + x_3 j^{2} a_2 + x_3^{2} (j a_2^{2} + j^{2} a_3^{2}) \ ,\  x_1 a_3  + x_2 a_3 + x_2^{2} (a_2^{2} + j a_3^{2})  + x_3 (a_1^{2} + a_2^{2} + j^{2} a_3+j a_3^{2})  + x_3^{2} (j^{2} a_2^{2} + j^{2} a_3^{2}) )		$
\item{\bf	XXXII: (2)	}	$		    (x_1 a_1  + x_2 (a_1+a_1^{2})  + x_2^{2} j a_3^{2} + x_3 (a_2^{2} + j^{2} a_3^{2})  + x_3^{2} a_3^{2}\ ,\  x_1 a_2  + x_2 (a_1^{2} + a_2+a_2^{2})  + x_2^{2} j a_3^{2} + x_3^{2} a_3^{2}\ ,\  x_1 a_3  + x_2 a_3 + x_2^{2} (a_2^{2} + a_3^{2})  + x_3 (a_1^{2} + a_2^{2})  + x_3^{2} (j^{2} a_2^{2} + j a_3^{2}) )		$
\item{\bf	XXXIII: (1)	}	$		    (x_1 a_1  + x_2 j a_2^{2} + x_2^{2} (j^{6} a_1+j^{6} a_1^{2} + j^{6} a_2^{2})  + x_2^{4} (j^{6} a_1+j^{6} a_1^{2} + j^{3} a_2^{2}) \ ,\  x_1 a_2  + x_2 (a_1^{2} + j^{5} a_2^{2})  + x_2^{2} (j^{6} a_2+j a_2^{2})  + x_2^{4} (j^{6} a_2+j^{5} a_2^{2}) )		$
\item{\bf	XXXIV: (2)	}	$		    (x_1 a_1  + x_2 (a_1+a_1^{2} + a_2 + j a_3+a_3^{2})  + x_2^{2} (a_2^{2} + j^{2} a_3^{2})  + x_3 (j^{2} a_2+a_2^{2} + a_3)  + x_3^{2} j a_2^{2}\ ,\  x_1 a_2  + x_2 (a_1^{2} + a_2^{2} + j a_3)  + x_2^{2} (a_2^{2} + j^{2} a_3^{2})  + x_3 (j^{2} a_2 + a_3)  + x_3^{2} j a_2^{2}\ ,\  x_1 a_3  + x_2 (a_2 + j^{2} a_3)  + x_3 (a_1^{2} + a_2+a_2^{2} + j a_3)  + x_3^{2} (j^{2} a_2^{2} + j^{2} a_3^{2}) )		$
\item{\bf	XXXV: (2)	}	$		    (x_1 a_1  + x_2 (j^{2} a_1+j^{2} a_1^{2} + j a_2^{2})  + x_2^{2} (a_2+j^{2} a_2^{2} + a_3^{2})  + x_3 (j a_2^{2} + j a_3+j a_3^{2})  + x_3^{2} (j^{2} a_1+j^{2} a_1^{2} + j a_2+j^{2} a_2^{2} + a_3) \ ,\  x_1 a_2  + x_2 (a_1^{2} + j^{2} a_2+j a_2^{2})  + x_2^{2} a_2^{2} + x_3 (j^{2} a_2^{2} + j^{2} a_3)  + x_3^{2} (a_1+a_1^{2} + j^{2} a_2+j a_2^{2}) \ ,\  x_1 a_3  + x_2 (j^{2} a_2 + j^{2} a_3+a_3^{2})  + x_2^{2} (j a_1+j a_1^{2} + a_2+j^{2} a_2^{2})  + x_3 (a_1^{2} + j a_2^{2})  + x_3^{2} a_3^{2})		$
\item{\bf	XXXVI: (2)	}	$		    (x_1 a_1  + x_2 (j a_1+j a_1^{2} + a_2+j a_2^{2})  + x_2^{2} (j^{2} a_2 + a_3)  + x_3 (j^{2} a_1+j^{2} a_1^{2} + j^{2} a_3+j a_3^{2})  + x_3^{2} j a_2\ ,\  x_1 a_2  + x_2 (a_1^{2} + a_2+a_2^{2} + a_3+j^{2} a_3^{2})  + x_2^{2} j^{2} a_3 + x_3 (j a_2+j^{2} a_2^{2})  + x_3^{2} (a_2 + j a_3) \ ,\  x_1 a_3  + x_2 (a_3+j^{2} a_3^{2})  + x_2^{2} a_2 + x_3 (a_1^{2} + j^{2} a_2+a_2^{2} + j a_3+j^{2} a_3^{2})  + x_3^{2} a_3)		$
\item{\bf	XXXVII: (2)	}	$		    (x_1 a_1  + x_2 (j^{2} a_1+j^{2} a_1^{2} + a_2^{2} + j^{2} a_3^{2})  + x_2^{2} (j a_1+j a_1^{2} + j a_2+j^{2} a_2^{2} + j^{2} a_3+j^{2} a_3^{2})  + x_3 (a_1+a_1^{2} + j a_2^{2} + j a_3^{2})  + x_3^{2} (j a_1+j a_1^{2} + j a_2^{2} + a_3+j a_3^{2}) \ ,\  x_1 a_2  + x_2 (a_1^{2} + j^{2} a_2+j a_2^{2} + j a_3^{2})  + x_2^{2} (j^{2} a_1+j^{2} a_1^{2} + j^{2} a_2+j a_2^{2} + j a_3+j^{2} a_3^{2})  + x_3 (a_2+j^{2} a_2^{2})  + x_3^{2} (j a_1+j a_1^{2} + j a_2^{2} + a_3+j a_3^{2}) \ ,\  x_1 a_3  + x_2 (a_2^{2} + j^{2} a_3+j^{2} a_3^{2})  + x_2^{2} (j^{2} a_1+j^{2} a_1^{2} + j^{2} a_2+j^{2} a_2^{2} + a_3+j a_3^{2})  + x_3 (a_1^{2} + a_2^{2} + a_3+j a_3^{2})  + x_3^{2} (j^{2} a_1+j^{2} a_1^{2} + j^{2} a_2^{2} + j a_3+j^{2} a_3^{2}) )		$
\item{\bf	XXXVIII: (2)	}	$		    (x_1 a_1  + x_2 (a_1+a_1^{2} + j a_2^{2} + j^{2} a_3+j^{2} a_3^{2})  + x_2^{2} (j a_1+j a_1^{2} + j a_2^{2} + j^{2} a_3^{2})  + x_3 (a_1+a_1^{2} + a_3)  + x_3^{2} (a_1+a_1^{2} + j^{2} a_2^{2} + j^{2} a_3+a_3^{2}) \ ,\  x_1 a_2  + x_2 (a_1^{2} + a_2+j^{2} a_2^{2})  + x_2^{2} (j a_2 + j^{2} a_3^{2})  + x_3 (a_2+j^{2} a_2^{2})  + x_3^{2} (a_2 + j a_3) \ ,\  x_1 a_3  + x_2 (a_2^{2} + a_3+j^{2} a_3^{2})  + x_2^{2} (j a_3+j a_3^{2})  + x_3 (a_1^{2} + j^{2} a_2^{2} + a_3+a_3^{2})  + x_3^{2} a_3^{2})		$
\item{\bf	XXXIX: (2)	}	$		    (x_1 a_1  + x_2 (j a_1+j a_1^{2} + j a_2^{2} + j^{2} a_3+j a_3^{2})  + x_2^{2} (a_2+a_2^{2} + a_3+a_3^{2})  + x_3 (j a_2^{2} + j^{2} a_3+j a_3^{2})  + x_3^{2} (j a_2+j a_2^{2} + a_3) \ ,\  x_1 a_2  + x_2 (j a_1+j^{2} a_1^{2} + a_2^{2} + j^{2} a_3+j^{2} a_3^{2})  + x_2^{2} (j^{2} a_2+j^{2} a_2^{2})  + x_3 (j^{2} a_1+j^{2} a_1^{2} + j a_2 + a_3)  + x_3^{2} (j^{2} a_2 + a_3+a_3^{2}) \ ,\  x_1 a_3  + x_2 (a_1+a_1^{2} + a_2^{2} + j a_3+a_3^{2})  + x_2^{2} (j^{2} a_2+j a_2^{2} + j^{2} a_3+j^{2} a_3^{2})  + x_3 (a_1 + a_2^{2} + a_3^{2})  + x_3^{2} (a_2+j^{2} a_2^{2} + j^{2} a_3+j a_3^{2}) )		$
\item{\bf	XL: (2)	}	$		    (x_1 a_1  + x_2 (j^{2} a_2 + j^{2} a_3)  + x_2^{2} (j a_1+j a_1^{2} + j^{2} a_2^{2} + j a_3^{2})  + x_3 (j^{2} a_1+j^{2} a_1^{2} + j a_2+a_2^{2} + a_3+j^{2} a_3^{2})  + x_3^{2} (j^{2} a_1+j^{2} a_1^{2} + j a_2+j^{2} a_2^{2} + j a_3+j^{2} a_3^{2}) \ ,\  x_1 a_2  + x_2 (a_1 + j a_2)  + x_2^{2} (j a_2 + j^{2} a_3^{2})  + x_3 (a_3+a_3^{2})  + x_3^{2} (j^{2} a_1+j^{2} a_1^{2} + j a_2+j^{2} a_2^{2} + a_3^{2}) \ ,\  x_1 a_3  + x_2 j a_3 + x_2^{2} (a_2^{2} + j a_3)  + x_3 (a_1 + j^{2} a_2+j a_2^{2} + j^{2} a_3)  + x_3^{2} (j a_3+j^{2} a_3^{2}) )		$
\item{\bf	XLI: (2)	}	$		    (x_1 a_1  + x_2 (j^{2} a_2+j^{2} a_2^{2} + a_3+a_3^{2})  + x_3 (a_1+a_1^{2} + j a_2 + a_3+j a_3^{2}) \ ,\  x_1 a_2  + x_2 (a_1^{2} + a_2+a_2^{2} + j a_3)  + x_3 (j a_2^{2} + a_3^{2}) \ ,\  x_1 a_3  + x_2 (j a_1+j a_1^{2} + a_2 + j^{2} a_3+j a_3^{2})  + x_3 (a_1 + a_3+j^{2} a_3^{2}) )		$
\item{\bf	XLII: (2)	}	$		    (x_1 a_1  + x_2 (j^{2} a_1+j^{2} a_1^{2} + j a_2+a_2^{2} + j^{2} a_3+a_3^{2})  + x_2^{2} (j^{2} a_1+j^{2} a_1^{2} + j^{2} a_2^{2} + j^{2} a_3^{2})  + x_3 (j a_1+j a_1^{2} + j a_3+j a_3^{2})  + x_3^{2} (j a_1+j a_1^{2} + a_2+j^{2} a_2^{2} + a_3^{2}) \ ,\  x_1 a_2  + x_2 (j a_1+j^{2} a_1^{2} + j^{2} a_2+a_2^{2} + a_3^{2})  + x_2^{2} (j a_1+j a_1^{2} + a_2+j^{2} a_2^{2})  + x_3 (j a_2 + j^{2} a_3+j^{2} a_3^{2})  + x_3^{2} (a_1+a_1^{2} + j a_2+j^{2} a_2^{2} + j^{2} a_3^{2}) \ ,\  x_1 a_3  + x_2 (a_1+a_1^{2} + a_2 + j^{2} a_3+j^{2} a_3^{2})  + x_2^{2} j^{2} a_3 + x_3 (a_1^{2} + j^{2} a_3+j^{2} a_3^{2})  + x_3^{2} (j^{2} a_1+j^{2} a_1^{2} + j^{2} a_2 + j^{2} a_3+j a_3^{2}) )		$
\item{\bf	XLIII: (2)	}	$		    (x_1 a_1  + x_2 (j^{2} a_1+j^{2} a_1^{2} + a_2+a_2^{2} + j a_3+j a_3^{2})  + x_3 (a_1+a_1^{2} + j a_2^{2} + j a_3+a_3^{2}) \ ,\  x_1 a_2  + x_2 (a_1 + a_2+a_2^{2} + j a_3+a_3^{2})  + x_3 (j a_1+j a_1^{2} + a_2^{2} + a_3) \ ,\  x_1 a_3  + x_2 (j^{2} a_1+j^{2} a_1^{2} + j a_2+j^{2} a_2^{2} + j^{2} a_3^{2})  + x_3 (j^{2} a_1+j a_1^{2} + a_2^{2} + a_3+j^{2} a_3^{2}) )		$
\item{\bf	XLIV: (2)	}	$		    (x_1 a_1  + x_2 (j^{2} a_1+j^{2} a_1^{2} + j^{2} a_2+a_2^{2} + j a_3+j^{2} a_3^{2})  + x_2^{2} (a_2 + j^{2} a_3)  + x_3 (j a_2 + a_3)  + x_3^{2} (a_1+a_1^{2} + a_2+j a_2^{2} + j^{2} a_3+a_3^{2}) \ ,\  x_1 a_2  + x_2 (a_1^{2} + a_2 + j a_3)  + x_2^{2} (j^{2} a_1+j^{2} a_1^{2} + a_2+j^{2} a_2^{2} + a_3+j^{2} a_3^{2})  + x_3 (j a_1+j a_1^{2} + j a_2^{2} + j a_3^{2})  + x_3^{2} (j a_2+j a_2^{2} + j^{2} a_3+j a_3^{2}) \ ,\  x_1 a_3  + x_2 (j^{2} a_1+j^{2} a_1^{2} + a_2+j a_2^{2} + a_3+j^{2} a_3^{2})  + x_2^{2} j^{2} a_2 + x_3 (a_1 + j^{2} a_3)  + x_3^{2} (a_1+a_1^{2} + j a_2+j^{2} a_2^{2} + j a_3+a_3^{2}) )		$
\item{\bf	XLV:	 (2)}	$		    (x_1 a_1  + x_2 (j a_1+j a_1^{2} + j^{2} a_2+j^{2} a_2^{2} + j^{2} a_3+j a_3^{2})  + x_3 (j^{2} a_2+a_2^{2} + a_3+a_3^{2}) \ ,\  x_1 a_2  + x_2 (a_1^{2} + a_2+a_2^{2} + j a_3)  + x_3 (a_2^{2} + j^{2} a_3) \ ,\  x_1 a_3  + x_2 (a_2 + a_3+a_3^{2})  + x_3 (a_1^{2} + a_2 + a_3+j a_3^{2}) )		$
\item{\bf	XLVI: (2)}	$		    (x_1 a_1  + x_2 (j a_1+j a_1^{2} + j^{2} a_2^{2})  + x_2^{2} (j^{2} a_1+j^{2} a_1^{2} + j^{2} a_2^{2} + j a_3^{2})  + x_3 (j a_1+j a_1^{2} + j^{2} a_2^{2})  + x_3^{2} (a_2^{2} + j^{2} a_3^{2}) \ ,\  x_1 a_2  + x_2 (a_1^{2} + j a_2 + j^{2} a_3^{2})  + x_2^{2} (j^{2} a_2+a_2^{2} + a_3^{2})  + x_3 (j a_2+j a_2^{2} + j a_3^{2})  + x_3^{2} a_2^{2}\ ,\  x_1 a_3  + x_2 (j a_3+a_3^{2})  + x_2^{2} (a_2^{2} + j^{2} a_3)  + x_3 (a_1^{2} + a_2^{2} + j a_3) )		$
\item{\bf	XLVII: (2)	}	$		    (x_1 a_1  + x_2 (a_2 + j^{2} a_3)  + x_2^{2} (a_2^{2} + a_3+j a_3^{2})  + x_3 (a_2^{2} + j^{2} a_3+j a_3^{2})  + x_3^{2} (j a_1+j a_1^{2} + j^{2} a_2 + j a_3+j a_3^{2}) \ ,\  x_1 a_2  + x_2 (a_1^{2} + j a_2+j a_2^{2} + j^{2} a_3+j^{2} a_3^{2})  + x_2^{2} (a_1+a_1^{2} + a_2+a_2^{2} + a_3+j^{2} a_3^{2})  + x_3 (a_1+a_1^{2} + j^{2} a_2 + j a_3+j^{2} a_3^{2})  + x_3^{2} (a_1+a_1^{2} + a_2^{2} + j a_3+a_3^{2}) \ ,\  x_1 a_3  + x_2 (j^{2} a_2+a_2^{2} + j a_3+j^{2} a_3^{2})  + x_2^{2} (a_1+a_1^{2} + j a_2+a_2^{2} + j a_3+a_3^{2})  + x_3 (a_1 + a_2+j^{2} a_2^{2} + a_3^{2})  + x_3^{2} (a_1+a_1^{2} + a_2+a_2^{2} + j^{2} a_3+j^{2} a_3^{2}) )		$
\item{\bf	XLVIII: (2)	}	$		    (x_1 a_1  + x_2 (a_2^{2} + j a_3^{2})  + x_2^{2} (a_1+a_1^{2} + j a_2+j^{2} a_2^{2} + a_3+j a_3^{2})  + x_3 (j a_2+a_2^{2} + a_3)  + x_3^{2} (j a_1+j a_1^{2} + j^{2} a_2+j^{2} a_2^{2} + j a_3) \ ,\  x_1 a_2  + x_2 (j^{2} a_1+j a_1^{2} + j a_2+j a_2^{2} + j a_3+a_3^{2})  + x_2^{2} (j^{2} a_2+j^{2} a_2^{2} + a_3+j a_3^{2})  + x_3 (j a_1+j a_1^{2} + j^{2} a_2+j^{2} a_2^{2} + j a_3^{2})  + x_3^{2} (j a_2 + a_3+j^{2} a_3^{2}) \ ,\  x_1 a_3  + x_2 (a_1+a_1^{2} + j a_2+a_2^{2} + j^{2} a_3^{2})  + x_2^{2} (j^{2} a_2+a_2^{2} + j^{2} a_3+j^{2} a_3^{2})  + x_3 (j^{2} a_1+j a_1^{2} + j a_2 + j a_3)  + x_3^{2} (j a_2+j a_2^{2} + a_3) )		$
\item{\bf	XLIX: (2)	}	$		    (x_1 a_1  + x_2 (j^{2} a_1+j^{2} a_1^{2} + j^{2} a_2+a_2^{2} + a_3+j a_3^{2})  + x_2^{2} (a_2+j^{2} a_2^{2} + j a_3^{2})  + x_3 (j^{2} a_1+j^{2} a_1^{2} + a_2 + j^{2} a_3+j a_3^{2})  + x_3^{2} (j a_1+j a_1^{2} + j a_3^{2}) \ ,\  x_1 a_2  + x_2 (j^{2} a_1+j a_1^{2} + a_2+j^{2} a_2^{2} + j^{2} a_3+j a_3^{2})  + x_2^{2} (j^{2} a_1+j^{2} a_1^{2} + j a_2 + a_3+a_3^{2})  + x_3 (a_1+a_1^{2})  + x_3^{2} (j a_1+j a_1^{2} + a_2+j a_2^{2} + j^{2} a_3) \ ,\  x_1 a_3  + x_2 (a_1+a_1^{2} + j^{2} a_2+j a_2^{2} + j^{2} a_3^{2})  + x_2^{2} (j a_1+j a_1^{2} + j^{2} a_2+j^{2} a_2^{2} + j a_3+j a_3^{2})  + x_3 (a_1^{2} + a_2+a_2^{2} + j^{2} a_3+a_3^{2})  + x_3^{2} (j a_1+j a_1^{2} + j a_2 + j a_3+j a_3^{2}) )		$
\item{\bf	L: (2)	}	$		    (x_1 a_1  + x_2 (a_2+j^{2} a_2^{2} + a_3^{2})  + x_2^{2} (j^{2} a_1+j^{2} a_1^{2} + j a_2)  + x_3 (j^{2} a_1+j^{2} a_1^{2} + a_2+j^{2} a_2^{2} + j^{2} a_3^{2})  + x_3^{2} (a_1+a_1^{2} + j^{2} a_3+j^{2} a_3^{2}) \ ,\  x_1 a_2  + x_2 (a_1^{2} + a_2)  + x_2^{2} (j a_1+j a_1^{2} + a_2 + j^{2} a_3+a_3^{2})  + x_3 (j^{2} a_1+j^{2} a_1^{2} + j a_2+j a_2^{2} + a_3+j^{2} a_3^{2})  + x_3^{2} (a_2^{2} + a_3+j a_3^{2}) \ ,\  x_1 a_3  + x_2 (a_1+a_1^{2} + a_2 + a_3+j a_3^{2})  + x_2^{2} (a_3+j^{2} a_3^{2})  + x_3 (j a_1+j^{2} a_1^{2} + j a_2+a_2^{2} + j a_3+a_3^{2})  + x_3^{2} (j^{2} a_1+j^{2} a_1^{2} + j^{2} a_2 + j^{2} a_3+j^{2} a_3^{2}) )		$
\item{\bf	LI: (2)	}	$		    (x_1 a_1  + x_2 (j a_1+j a_1^{2} + a_3+j a_3^{2})  + x_2^{2} (a_2+a_2^{2} + a_3)  + x_3 (j a_1+j a_1^{2} + j a_2+a_2^{2} + j^{2} a_3^{2})  + x_3^{2} (j a_2+j a_2^{2} + j^{2} a_3^{2}) \ ,\  x_1 a_2  + x_2 (a_1^{2} + j a_2+j^{2} a_2^{2} + j^{2} a_3^{2})  + x_2^{2} (j^{2} a_2+j a_2^{2} + j^{2} a_3^{2})  + x_3 j a_2 + x_3^{2} (j^{2} a_2+j a_2^{2} + j a_3) \ ,\  x_1 a_3  + x_2 j a_3 + x_2^{2} (a_2^{2} + j^{2} a_3+j^{2} a_3^{2})  + x_3 (a_1^{2} + j^{2} a_2^{2} + j a_3+j^{2} a_3^{2})  + x_3^{2} (j a_2+j^{2} a_2^{2} + j a_3+j a_3^{2}) )		$
\item{\bf	LII: (2)	}	$		    (x_1 a_1  + x_2 (a_1+a_1^{2} + j^{2} a_2+j a_2^{2} + j a_3)  + x_2^{2} (j a_1+j a_1^{2} + j^{2} a_2+j a_2^{2})  + x_3 (j^{2} a_1+j^{2} a_1^{2} + j a_3+a_3^{2})  + x_3^{2} (j a_1+j a_1^{2} + j a_2^{2}) \ ,\  x_1 a_2  + x_2 (j^{2} a_1+j a_1^{2} + j^{2} a_2 + j a_3)  + x_2^{2} (a_1+a_1^{2} + a_2+j a_2^{2} + a_3)  + x_3 (a_1+a_1^{2} + a_2+j a_2^{2} + j a_3+a_3^{2})  + x_3^{2} (j a_1+j a_1^{2} + j a_2^{2} + j^{2} a_3+a_3^{2}) \ ,\  x_1 a_3  + x_2 (j^{2} a_1+j^{2} a_1^{2} + j^{2} a_2+j a_2^{2} + j a_3+a_3^{2})  + x_2^{2} (j a_1+j a_1^{2} + j a_2+j a_2^{2} + a_3^{2})  + x_3 (a_1 + j a_2+j a_2^{2} + a_3+a_3^{2})  + x_3^{2} (j a_2+j a_2^{2} + j a_3^{2}) )		$
\item{\bf	LIII: (2)	}	$		    (x_1 a_1  + x_2 (j^{2} a_1+j^{2} a_1^{2} + j^{2} a_2+j^{2} a_2^{2} + a_3)  + x_3 (j a_1+j a_1^{2} + j a_2 + j a_3+j a_3^{2}) \ ,\  x_1 a_2  + x_2 (j a_1+j^{2} a_1^{2} + j^{2} a_2+j^{2} a_2^{2} + j a_3)  + x_3 (j^{2} a_1+j^{2} a_1^{2} + a_2^{2} + j^{2} a_3^{2}) \ ,\  x_1 a_3  + x_2 (a_1+a_1^{2} + a_2 + j a_3+a_3^{2})  + x_3 (a_1 + j^{2} a_2+j^{2} a_2^{2} + a_3^{2}) )		$
\item{\bf	LIV: (2)	}	$		    (x_1 a_1  + x_2 (j^{2} a_1+j^{2} a_1^{2} + a_3)  + x_3 (a_1+a_1^{2} + j a_2^{2} + j a_3+a_3^{2}) \ ,\  x_1 a_2  + x_2 (j a_1+j^{2} a_1^{2} + j^{2} a_2+j^{2} a_2^{2} + j^{2} a_3)  + x_3 (j a_1+j a_1^{2} + a_2^{2} + j^{2} a_3+j^{2} a_3^{2}) \ ,\  x_1 a_3  + x_2 (a_1+a_1^{2} + a_2 + j^{2} a_3+a_3^{2})  + x_3 (a_1^{2} + j^{2} a_2+j^{2} a_2^{2} + a_3+j a_3^{2}) )		$
\item{\bf	LV: (2)	}	$		    (x_1 a_1  + x_2 (j a_1+j a_1^{2} + j a_2+j^{2} a_2^{2} + j a_3^{2})  + x_2^{2} (a_1+a_1^{2} + j^{2} a_2+j a_2^{2} + a_3^{2})  + x_3 (a_1+a_1^{2} + j a_2^{2} + a_3^{2})  + x_3^{2} j a_2\ ,\  x_1 a_2  + x_2 (j a_1+j^{2} a_1^{2} + j a_2+a_2^{2} + j^{2} a_3^{2})  + x_2^{2} j^{2} a_2 + x_3 (j a_1+j a_1^{2} + j^{2} a_2^{2} + j a_3^{2})  + x_3^{2} j a_2\ ,\  x_1 a_3  + x_2 (j a_1+j a_1^{2} + a_2^{2} + a_3)  + x_2^{2} (j a_2+j a_2^{2} + a_3+a_3^{2})  + x_3 (j^{2} a_1+j a_1^{2} + j a_2+a_2^{2} + j^{2} a_3)  + x_3^{2} (j^{2} a_2+j^{2} a_2^{2}) )		$
\item{\bf	LVI: (2)	}	$		    (x_1 a_1  + x_2 (j a_2 + j^{2} a_3+j^{2} a_3^{2})  + x_2^{2} (j a_2^{2} + a_3)  + x_3 (a_1+a_1^{2} + a_2 + j^{2} a_3+j^{2} a_3^{2})  + x_3^{2} (j a_1+j a_1^{2} + j^{2} a_2+j^{2} a_2^{2} + a_3+j^{2} a_3^{2}) \ ,\  x_1 a_2  + x_2 (j^{2} a_1+j a_1^{2} + j a_2 + j a_3+j^{2} a_3^{2})  + x_2^{2} (j a_1+j a_1^{2} + j^{2} a_2+a_2^{2} + j^{2} a_3^{2})  + x_3 (j a_1+j a_1^{2} + a_2^{2} + j^{2} a_3+j^{2} a_3^{2})  + x_3^{2} (j a_2+j^{2} a_2^{2} + j^{2} a_3+j a_3^{2}) \ ,\  x_1 a_3  + x_2 (j^{2} a_2+j a_2^{2} + a_3)  + x_2^{2} (a_1+a_1^{2} + j a_2+j a_2^{2} + j^{2} a_3+a_3^{2})  + x_3 (j^{2} a_1+j a_1^{2} + j^{2} a_2+a_2^{2} + j^{2} a_3+a_3^{2})  + x_3^{2} (j^{2} a_2+j a_2^{2} + j a_3) )		$
\item{\bf	LVII:	(2)}	$		    (x_1 a_1  + x_2 (j^{2} a_1+j^{2} a_1^{2} + j^{2} a_2+j a_2^{2} + j a_3+j^{2} a_3^{2})  + x_2^{2} (j a_1+j a_1^{2} + a_2^{2} + a_3^{2})  + x_3 (j^{2} a_2+j a_2^{2} + j a_3+j^{2} a_3^{2})  + x_3^{2} (j a_1+j a_1^{2} + j^{2} a_2^{2} + a_3+j^{2} a_3^{2}) \ ,\  x_1 a_2  + x_2 (a_1 + j^{2} a_2+a_2^{2} + j^{2} a_3+j^{2} a_3^{2})  + x_2^{2} (a_1+a_1^{2} + j^{2} a_2+a_2^{2} + a_3^{2})  + x_3 (a_2+a_2^{2} + j^{2} a_3+j^{2} a_3^{2})  + x_3^{2} (j^{2} a_1+j^{2} a_1^{2} + j^{2} a_2+j a_2^{2} + j a_3^{2}) \ ,\  x_1 a_3  + x_2 (j^{2} a_1+j^{2} a_1^{2} + j a_2+j^{2} a_2^{2})  + x_2^{2} (j a_2+j a_2^{2} + j^{2} a_3)  + x_3 (j^{2} a_1+j a_1^{2} + j a_3)  + x_3^{2} (j^{2} a_2 + j a_3+j^{2} a_3^{2}) )		$
\item{\bf	LVIII: (2)	}	$		    (x_1 a_1  + x_2 (j a_2+j a_2^{2} + j^{2} a_3+j a_3^{2})  + x_3 (a_1+a_1^{2} + a_2 + j^{2} a_3+a_3^{2}) \ ,\  x_1 a_2  + x_2 (a_1 + j a_2+j a_2^{2} + j a_3^{2})  + x_3 (a_1+a_1^{2} + j a_2+j a_2^{2} + j a_3+j a_3^{2}) \ ,\  x_1 a_3  + x_2 (j^{2} a_1+j^{2} a_1^{2} + a_2 + j^{2} a_3^{2})  + x_3 (j a_1+j^{2} a_1^{2} + a_2+j^{2} a_2^{2} + j a_3+j a_3^{2}) )		$
\item{\bf	LIX:	(2)}	$		    (x_1 a_1  + x_2 (j^{2} a_1+j^{2} a_1^{2} + j^{2} a_2+a_2^{2} + a_3+j a_3^{2})  + x_2^{2} (j^{2} a_1+j^{2} a_1^{2} + j a_2^{2} + j^{2} a_3^{2})  + x_3 (j^{2} a_1+j^{2} a_1^{2} + j a_2+j^{2} a_2^{2} + j^{2} a_3+j a_3^{2})  + x_3^{2} (j a_2+j^{2} a_2^{2} + j^{2} a_3) \ ,\  x_1 a_2  + x_2 (j^{2} a_1+j a_1^{2} + j^{2} a_2^{2} + j a_3+j^{2} a_3^{2})  + x_2^{2} (j a_1+j a_1^{2} + j^{2} a_2^{2} + j^{2} a_3)  + x_3 (j a_1+j a_1^{2} + j^{2} a_2)  + x_3^{2} (j^{2} a_1+j^{2} a_1^{2} + j^{2} a_2^{2} + a_3) \ ,\  x_1 a_3  + x_2 (j^{2} a_1+j^{2} a_1^{2} + a_2+j a_2^{2} + j^{2} a_3+a_3^{2})  + x_2^{2} (j^{2} a_1+j^{2} a_1^{2} + a_2+j^{2} a_2^{2})  + x_3 (a_1^{2} + a_2+j a_2^{2} + a_3+j a_3^{2})  + x_3^{2} (j a_1+j a_1^{2} + j^{2} a_2+j a_2^{2} + j a_3+j a_3^{2}) )		$
\item{\bf	LX: (2)	}	$		    (x_1 a_1  + x_2 (a_1+a_1^{2} + j a_2+j^{2} a_2^{2} + j^{2} a_3+a_3^{2})  + x_2^{2} (j a_1+j a_1^{2} + j a_2+j^{2} a_2^{2} + a_3+a_3^{2})  + x_3 (a_1+a_1^{2} + j^{2} a_2+a_2^{2} + j^{2} a_3)  + x_3^{2} (j^{2} a_1+j^{2} a_1^{2} + j a_2^{2} + j a_3+a_3^{2}) \ ,\  x_1 a_2  + x_2 (j a_1+j^{2} a_1^{2} + a_2^{2} + j a_3)  + x_2^{2} (j^{2} a_1+j^{2} a_1^{2} + a_2 + j^{2} a_3+j^{2} a_3^{2})  + x_3 (j^{2} a_1+j^{2} a_1^{2} + j^{2} a_2^{2} + j a_3+a_3^{2})  + x_3^{2} (j^{2} a_1+j^{2} a_1^{2} + j a_2+j a_2^{2} + a_3^{2}) \ ,\  x_1 a_3  + x_2 (j^{2} a_1+j^{2} a_1^{2} + j^{2} a_2+a_2^{2} + j^{2} a_3+j a_3^{2})  + x_2^{2} (j^{2} a_1+j^{2} a_1^{2} + j a_2+a_2^{2} + j a_3+a_3^{2})  + x_3 (j a_1+j^{2} a_1^{2} + j a_2+j^{2} a_2^{2} + j a_3+j^{2} a_3^{2})  + x_3^{2} (j^{2} a_1+j^{2} a_1^{2} + j^{2} a_2 + a_3) )		$
\item{\bf	LXI: (2)	}	$		    (x_1 a_1  + x_2 (a_2+j a_2^{2} + j^{2} a_3+j a_3^{2})  + x_2^{2} (j a_1+j a_1^{2} + j^{2} a_2 + a_3+j a_3^{2})  + x_3 (j^{2} a_1+j^{2} a_1^{2} + j^{2} a_2+j^{2} a_2^{2} + j^{2} a_3)  + x_3^{2} (j^{2} a_1+j^{2} a_1^{2} + j a_2+a_2^{2} + j a_3^{2}) \ ,\  x_1 a_2  + x_2 (a_1^{2} + j a_2+j a_2^{2} + j a_3^{2})  + x_2^{2} (a_1+a_1^{2} + j a_2+j a_2^{2} + j^{2} a_3)  + x_3 (a_2 + j a_3+j^{2} a_3^{2})  + x_3^{2} (a_1+a_1^{2} + j^{2} a_2 + j^{2} a_3+j^{2} a_3^{2}) \ ,\  x_1 a_3  + x_2 (j a_2+j a_2^{2} + a_3^{2})  + x_2^{2} (j^{2} a_1+j^{2} a_1^{2} + j^{2} a_2+j a_2^{2} + a_3^{2})  + x_3 (j^{2} a_1+j a_1^{2} + a_2 + a_3+j^{2} a_3^{2})  + x_3^{2} (a_3+j a_3^{2}) )		$
\item{\bf	LXII:	(2)}	$		    (x_1 a_1  + x_2 (j a_1+j a_1^{2} + a_2^{2} + a_3^{2})  + x_2^{2} (j a_1+j a_1^{2} + j a_2+j^{2} a_2^{2} + j a_3+a_3^{2})  + x_3 (a_1+a_1^{2} + a_2 + a_3+a_3^{2})  + x_3^{2} (j a_1+j a_1^{2} + j^{2} a_2 + j^{2} a_3+a_3^{2}) \ ,\  x_1 a_2  + x_2 (a_1^{2} + j a_2+j a_2^{2})  + x_2^{2} (a_1+a_1^{2} + a_2+a_2^{2} + a_3)  + x_3 (j^{2} a_1+j^{2} a_1^{2} + j^{2} a_2+a_2^{2} + j^{2} a_3+j^{2} a_3^{2})  + x_3^{2} (j a_1+j a_1^{2} + j^{2} a_2+j a_2^{2} + j a_3+a_3^{2}) \ ,\  x_1 a_3  + x_2 (a_2^{2} + j a_3+j^{2} a_3^{2})  + x_2^{2} (j a_1+j a_1^{2} + j a_2+j a_2^{2} + j a_3+j a_3^{2})  + x_3 (a_1 + a_2+j^{2} a_2^{2} + a_3+j a_3^{2})  + x_3^{2} (j^{2} a_1+j^{2} a_1^{2} + j^{2} a_2+j^{2} a_2^{2} + j a_3+a_3^{2}) )		$
\item{\bf	LXIII: (2)	}	$		    (x_1 a_1  + x_2 (a_1+a_1^{2} + a_2+a_2^{2} + j a_3+j a_3^{2})  + x_3 (j a_2+a_2^{2} + j a_3+j^{2} a_3^{2}) \ ,\  x_1 a_2  + x_2 (a_1 + a_2+a_2^{2} + j a_3)  + x_3 (j a_1+j a_1^{2} + j^{2} a_2 + a_3^{2}) \ ,\  x_1 a_3  + x_2 (a_1+a_1^{2} + a_2^{2} + j a_3+j^{2} a_3^{2})  + x_3 (j^{2} a_1+j a_1^{2} + j^{2} a_2 + j^{2} a_3) )		$
\item{\bf	LXIV: (2)	}	$		    (x_1 a_1  + x_2 (j a_1+j a_1^{2} + j^{2} a_2 + j^{2} a_3+j a_3^{2})  + x_2^{2} (a_1+a_1^{2} + j^{2} a_2^{2} + j^{2} a_3^{2})  + x_3 (a_1+a_1^{2} + a_2 + a_3+j^{2} a_3^{2})  + x_3^{2} (a_1+a_1^{2} + j^{2} a_2+j a_2^{2} + j^{2} a_3) \ ,\  x_1 a_2  + x_2 (a_1^{2} + j a_2^{2})  + x_2^{2} (j^{2} a_1+j^{2} a_1^{2} + j a_2 + j a_3+a_3^{2})  + x_3 (a_1+a_1^{2} + j a_2^{2})  + x_3^{2} (j^{2} a_2+j^{2} a_2^{2} + j a_3) \ ,\  x_1 a_3  + x_2 (a_2+j^{2} a_2^{2} + a_3+a_3^{2})  + x_2^{2} (j^{2} a_1+j^{2} a_1^{2} + j^{2} a_2 + j^{2} a_3+a_3^{2})  + x_3 (j^{2} a_1+j a_1^{2} + j a_2 + j a_3+j^{2} a_3^{2})  + x_3^{2} (a_2+a_2^{2} + j a_3^{2}) )		$
\item{\bf	LXV:	 (2)}	$		    (x_1 a_1  + x_2 (a_2+j^{2} a_2^{2} + j a_3+j a_3^{2})  + x_2^{2} (j^{2} a_2+a_2^{2} + j a_3^{2})  + x_3 (a_1+a_1^{2} + a_2 + j^{2} a_3^{2})  + x_3^{2} (j^{2} a_1+j^{2} a_1^{2} + j a_2+j^{2} a_2^{2} + a_3) \ ,\  x_1 a_2  + x_2 (a_1 + a_2+j a_2^{2} + a_3+a_3^{2})  + x_2^{2} (j a_2+a_2^{2})  + x_3 (j^{2} a_1+j^{2} a_1^{2} + j^{2} a_2+j^{2} a_2^{2})  + x_3^{2} (j^{2} a_1+j^{2} a_1^{2} + j^{2} a_2+j a_2^{2} + j^{2} a_3) \ ,\  x_1 a_3  + x_2 (j^{2} a_1+j^{2} a_1^{2} + j a_2^{2} + j a_3)  + x_2^{2} (j^{2} a_1+j^{2} a_1^{2} + a_2+j a_2^{2} + j a_3+a_3^{2})  + x_3 (j a_1+j^{2} a_1^{2} + j a_2^{2} + j a_3^{2})  + x_3^{2} (j^{2} a_1+j^{2} a_1^{2} + j a_2+j^{2} a_2^{2} + a_3+a_3^{2}) )		$
\item{\bf	LXVI: (2)	}	$		    (x_1 a_1  + x_2 (j a_2^{2} + j a_3+a_3^{2})  + x_2^{2} (j a_1+j a_1^{2} + j^{2} a_2+a_2^{2} + j^{2} a_3)  + x_3 (a_1+a_1^{2} + a_2+j^{2} a_2^{2} + a_3)  + x_3^{2} (j a_1+j a_1^{2} + j a_2+j a_2^{2} + a_3+j a_3^{2}) \ ,\  x_1 a_2  + x_2 (a_1^{2} + j a_2^{2} + a_3)  + x_2^{2} (j a_2 + j a_3+a_3^{2})  + x_3 (j a_1+j a_1^{2} + j^{2} a_3+j^{2} a_3^{2})  + x_3^{2} (j a_1+j a_1^{2} + a_2^{2} + j^{2} a_3+a_3^{2}) \ ,\  x_1 a_3  + x_2 (j^{2} a_1+j^{2} a_1^{2} + j a_2^{2} + j^{2} a_3+a_3^{2})  + x_2^{2} (j^{2} a_1+j^{2} a_1^{2} + a_2+j a_2^{2} + a_3+a_3^{2})  + x_3 (j a_1+j^{2} a_1^{2} + j a_2+j^{2} a_2^{2} + a_3+j^{2} a_3^{2})  + x_3^{2} (j^{2} a_1+j^{2} a_1^{2} + j^{2} a_2 + j^{2} a_3+j^{2} a_3^{2}) )		$
\item{\bf	LXVII: (2)	}	$		    (x_1 a_1  + x_2 (j a_1+j a_1^{2} + j a_2+j a_2^{2} + a_3^{2})  + x_3 (a_1+a_1^{2} + j^{2} a_2+j a_2^{2} + j^{2} a_3+j^{2} a_3^{2}) \ ,\  x_1 a_2  + x_2 (a_1 + j^{2} a_2+j^{2} a_2^{2} + j a_3+j^{2} a_3^{2})  + x_3 (j^{2} a_1+j^{2} a_1^{2} + j a_2 + j^{2} a_3^{2}) \ ,\  x_1 a_3  + x_2 (a_2 + j a_3)  + x_3 (a_1^{2} + a_2+j a_2^{2} + j a_3+j^{2} a_3^{2}) )		$
\item{\bf	LXVIII: (2)	}	$		    (x_1 a_1  + x_2 (j^{2} a_1+j^{2} a_1^{2} + j a_2+j a_2^{2} + a_3+a_3^{2})  + x_3 (a_1+a_1^{2} + a_2+j^{2} a_2^{2} + j a_3+a_3^{2}) \ ,\  x_1 a_2  + x_2 (j a_1+j^{2} a_1^{2} + a_2+a_2^{2} + j a_3^{2})  + x_3 (j^{2} a_1+j^{2} a_1^{2} + a_2^{2} + j a_3^{2}) \ ,\  x_1 a_3  + x_2 (a_1+a_1^{2} + a_2^{2} + j a_3^{2})  + x_3 (a_1 + a_2+a_2^{2}) )		$
\item{\bf	LXIX: (2)	}	$		    (x_1 a_1  + x_2 (j^{2} a_1+j^{2} a_1^{2} + a_2+a_2^{2} + j a_3+j a_3^{2})  + x_3 (a_1+a_1^{2} + j^{2} a_2+a_2^{2} + j a_3+a_3^{2}) \ ,\  x_1 a_2  + x_2 (a_1^{2} + a_2+a_2^{2} + j^{2} a_3)  + x_3 (a_1+a_1^{2} + j^{2} a_2+j^{2} a_2^{2} + a_3+j a_3^{2}) \ ,\  x_1 a_3  + x_2 (j a_1+j a_1^{2} + j a_2+j^{2} a_2^{2} + j^{2} a_3^{2})  + x_3 (j a_1+j^{2} a_1^{2} + a_2^{2} + j^{2} a_3+j a_3^{2}) )		$
\item{\bf	LXX:	 (2)}	$		    (x_1 a_1  + x_2 (j a_1+j a_1^{2} + j^{2} a_2 + j^{2} a_3^{2})  + x_2^{2} (j a_1+j a_1^{2} + j^{2} a_2 + a_3^{2})  + x_3 (a_1+a_1^{2} + j^{2} a_2)  + x_3^{2} (j a_1+j a_1^{2} + a_2 + j a_3^{2}) \ ,\  x_1 a_2  + x_2 (j^{2} a_1+j a_1^{2} + j^{2} a_2+j a_2^{2} + j^{2} a_3+j^{2} a_3^{2})  + x_2^{2} (a_1+a_1^{2} + a_2 + a_3+j^{2} a_3^{2})  + x_3 (j a_1+j a_1^{2} + a_2+j^{2} a_2^{2} + j a_3+j a_3^{2})  + x_3^{2} (a_2+j^{2} a_2^{2} + j^{2} a_3^{2}) \ ,\  x_1 a_3  + x_2 (j a_2^{2} + j a_3+j^{2} a_3^{2})  + x_2^{2} (j^{2} a_1+j^{2} a_1^{2} + a_2+j a_2^{2} + j a_3+j^{2} a_3^{2})  + x_3 (j a_1+j^{2} a_1^{2} + j^{2} a_2+a_2^{2} + a_3+j a_3^{2})  + x_3^{2} (j a_1+j a_1^{2} + j^{2} a_2+j a_2^{2} + j a_3+j^{2} a_3^{2}) )		$
\item{\bf	LXXI: (3)}	$			    (x_1 (j^{5} a_1 + j^{24} a_1^{2} + j^{43} a_1^{4} + j^{18} a_1^{8} + j^{58} a_1^{16} + j^{43} a_1^{32}) + x_1^{2} (j^{2} a_1 + j^{62} a_1^{2} + j^{45} a_1^{4} + j^{53} a_1^{8} + j^{14} a_1^{16} + j a_1^{32}) + x_1^{4} (j^{16} a_1 + j^{16} a_1^{2} + j^{45} a_1^{4} + j^{49} a_1^{8} + j^{50} a_1^{16} + j^{43} a_1^{32}) + x_1^{8} (j^{50} a_1 + j^{42} a_1^{4} + j^{4} a_1^{8} + j^{57} a_1^{16} + j^{5} a_1^{32}) + x_1^{16} (j^{40} a_1 + j^{16} a_1^{2} + j^{38} a_1^{4} + j^{18} a_1^{8} + j^{17} a_1^{16} + j^{27} a_1^{32}) + x_1^{32} (j^{44} a_1 + j^{10} a_1^{2} + j^{36} a_1^{4} + j^{5} a_1^{8} + j^{26} a_1^{16} + j^{37} a_1^{32}))		$
\item{\bf	LXXII: (2)	}	$		    (x_1 a_1  + x_2 (j a_2+a_2^{2} + j a_3+j^{2} a_3^{2})  + x_2^{2} (a_2+j a_2^{2} + j a_3+a_3^{2})  + x_3 (a_1+a_1^{2} + j^{2} a_2+j^{2} a_2^{2} + j a_3+a_3^{2})  + x_3^{2} (j a_2+a_2^{2} + j a_3+j a_3^{2}) \ ,\  x_1 a_2  + x_2 (a_1 + j a_2+j^{2} a_2^{2} + a_3+j a_3^{2})  + x_2^{2} (a_2^{2} + a_3)  + x_3 (j a_2 + a_3)  + x_3^{2} (a_1+a_1^{2} + j^{2} a_2^{2} + j a_3+j^{2} a_3^{2}) \ ,\  x_1 a_3  + x_2 (a_2 + j a_3)  + x_3 (a_1^{2} + j a_2+j^{2} a_2^{2} + a_3+a_3^{2})  + x_3^{2} (a_2 + j a_3+a_3^{2}) )		$
\item{\bf	LXXIII: (3)	}	$			    (x_1 (j^{57} a_1 + j^{16} a_1^{2} + j^{50} a_1^{4} + j^{37} a_1^{8} + j^{39} a_1^{16} + j a_1^{32}) + x_1^{2} (j^{29} a_1 + j^{38} a_1^{2} + j^{8} a_1^{4} + j^{30} a_1^{8} + j^{35} a_1^{16} + j^{54} a_1^{32}) + x_1^{4} (j^{15} a_1 + j^{4} a_1^{2} + j^{26} a_1^{4} + j^{32} a_1^{8} + j^{39} a_1^{32}) + x_1^{8} (j^{56} a_1 + j^{37} a_1^{2} + j^{23} a_1^{4} + j^{35} a_1^{8} + j^{6} a_1^{16} + j^{33} a_1^{32}) + x_1^{16} (j^{4} a_1^{2} + j^{55} a_1^{4} + j^{4} a_1^{8} + j^{4} a_1^{16} + j^{57} a_1^{32}) + x_1^{32} (j^{10} a_1 + j^{34} a_1^{2} + j^{29} a_1^{4} + j^{48} a_1^{8} + j^{55} a_1^{16} + j^{18} a_1^{32}))		$
\item{\bf	LXXIV: (3)	}	$			    (x_1 (j^{21} a_1 + j^{46} a_1^{2} + j^{25} a_1^{4} + j^{53} a_1^{8} + j^{44} a_1^{16} + j^{44} a_1^{32}) + x_1^{2} (j^{28} a_1 + j^{26} a_1^{2} + j^{10} a_1^{4} + j^{61} a_1^{8} + j^{34} a_1^{16} + j^{57} a_1^{32}) + x_1^{4} (j^{10} a_1 + j^{18} a_1^{2} + j^{57} a_1^{4} + j^{48} a_1^{8} + j^{33} a_1^{16} + j^{19} a_1^{32}) + x_1^{8} (j^{23} a_1 + j^{44} a_1^{2} + j^{37} a_1^{4} + j^{7} a_1^{8} + j^{40} a_1^{16} + a_1^{32}) + x_1^{16} (j^{6} a_1 + j^{14} a_1^{2} + j^{48} a_1^{4} + j^{17} a_1^{8} + j^{57} a_1^{16} + j a_1^{32}) + x_1^{32} (j^{20} a_1 + j^{6} a_1^{2} + j^{48} a_1^{4} + j^{39} a_1^{8} + j^{13} a_1^{16} + j^{27} a_1^{32}))		$
\item{\bf	LXXV: (3)	}	$			    (x_1 (j^{53} a_1 + j^{33} a_1^{2} + j^{14} a_1^{4} + j^{18} a_1^{8} + j^{29} a_1^{16} + j^{49} a_1^{32}) + x_1^{2} (j^{22} a_1 + j^{30} a_1^{2} + j^{50} a_1^{4} + j^{53} a_1^{8} + j^{52} a_1^{16} + j^{48} a_1^{32}) + x_1^{4} (j^{31} a_1 + j^{2} a_1^{2} + j^{47} a_1^{4} + j^{48} a_1^{8} + j^{10} a_1^{16} + j^{21} a_1^{32}) + x_1^{8} (j^{18} a_1 + j^{25} a_1^{2} + j^{26} a_1^{4} + a_1^{8} + j^{30} a_1^{16} + j^{9} a_1^{32}) + x_1^{16} (j^{16} a_1 + j^{52} a_1^{2} + j^{29} a_1^{4} + j^{10} a_1^{8} + j^{38} a_1^{16} + j^{52} a_1^{32}) + x_1^{32} (j^{47} a_1 + j^{33} a_1^{2} + j^{29} a_1^{4} + j^{45} a_1^{8} + j^{22} a_1^{16} + j^{62} a_1^{32}))		$
\item{\bf	LXXVI: (2)	}	$		    (x_1 a_1  + x_2 (j^{2} a_1+j^{2} a_1^{2} + j a_2+a_2^{2} + a_3+j^{2} a_3^{2})  + x_2^{2} (a_1+a_1^{2} + a_2+j a_2^{2} + j^{2} a_3)  + x_3 (a_1+a_1^{2} + j^{2} a_2^{2})  + x_3^{2} (j a_1+j a_1^{2} + j a_3^{2}) \ ,\  x_1 a_2  + x_2 (a_1^{2} + j^{2} a_2 + j^{2} a_3+j a_3^{2})  + x_2^{2} (a_2+j a_2^{2} + j a_3+a_3^{2})  + x_3 (a_2+j^{2} a_2^{2})  + x_3^{2} (j a_2+j^{2} a_2^{2} + a_3^{2}) \ ,\  x_1 a_3  + x_2 (j a_1+j a_1^{2} + j^{2} a_2+a_2^{2} + j^{2} a_3+a_3^{2})  + x_2^{2} (a_1+a_1^{2} + j a_2+a_2^{2} + a_3+j^{2} a_3^{2})  + x_3 (a_1^{2} + a_3)  + x_3^{2} (j a_2^{2} + j a_3+a_3^{2}) )		$
\item{\bf	LXXVII: (3)	}	$			    (x_1 (j^{53} a_1 + j^{38} a_1^{2} + j^{5} a_1^{4} + j^{39} a_1^{8} + j^{7} a_1^{16} + j^{5} a_1^{32}) + x_1^{2} (j^{14} a_1 + j^{26} a_1^{2} + j^{17} a_1^{4} + j^{27} a_1^{8} + j^{6} a_1^{16} + j^{29} a_1^{32}) + x_1^{4} (j^{38} a_1 + j^{44} a_1^{2} + j^{42} a_1^{4} + j^{48} a_1^{8} + j^{19} a_1^{16} + j^{56} a_1^{32}) + x_1^{8} (j^{30} a_1 + j^{7} a_1^{2} + j^{10} a_1^{4} + j^{26} a_1^{8} + j^{7} a_1^{16} + j^{2} a_1^{32}) + x_1^{16} (j^{15} a_1 + j^{28} a_1^{2} + j^{34} a_1^{4} + j^{61} a_1^{8} + j^{55} a_1^{16} + j^{18} a_1^{32}) + x_1^{32} (j^{58} a_1 + j^{51} a_1^{2} + j^{15} a_1^{4} + j^{6} a_1^{8} + j^{30} a_1^{16} + j^{8} a_1^{32}))		$
\item{\bf	LXXVIII: (3)	}	$			    (x_1 (j^{27} a_1 + j^{10} a_1^{2} + j^{39} a_1^{4} + j^{13} a_1^{8} + j^{22} a_1^{16} + j^{39} a_1^{32}) + x_1^{2} (j^{4} a_1 + j^{48} a_1^{2} + j^{14} a_1^{4} + j^{23} a_1^{8} + j^{47} a_1^{16} + j^{26} a_1^{32}) + x_1^{4} (j^{11} a_1 + j^{45} a_1^{2} + j^{35} a_1^{4} + j^{38} a_1^{8} + j^{30} a_1^{16} + j^{45} a_1^{32}) + x_1^{8} (j^{49} a_1 + j a_1^{2} + j^{11} a_1^{4} + j^{4} a_1^{8} + j^{58} a_1^{16} + j^{14} a_1^{32}) + x_1^{16} (j^{49} a_1 + j^{48} a_1^{2} + j^{34} a_1^{4} + j^{18} a_1^{8} + j^{46} a_1^{16} + j^{27} a_1^{32}) + x_1^{32} (j^{20} a_1 + j^{51} a_1^{2} + j^{55} a_1^{4} + j^{47} a_1^{8} + j^{37} a_1^{16} + j^{46} a_1^{32}))		$
\item{\bf	LXXIX: (3)	}	$			    (x_1 (j^{31} a_1 + j^{32} a_1^{2} + j^{29} a_1^{4} + j^{57} a_1^{8} + j^{26} a_1^{16} + j^{37} a_1^{32}) + x_1^{2} (j^{62} a_1 + j^{40} a_1^{2} + j^{55} a_1^{4} + j^{3} a_1^{8} + j^{15} a_1^{16}) + x_1^{4} (j^{31} a_1 + j^{3} a_1^{2} + j^{60} a_1^{4} + j^{29} a_1^{8} + j^{16} a_1^{16} + j^{58} a_1^{32}) + x_1^{8} (j^{9} a_1 + j^{52} a_1^{2} + j^{54} a_1^{4} + j^{5} a_1^{16} + j^{42} a_1^{32}) + x_1^{16} (j^{51} a_1 + j a_1^{2} + j^{43} a_1^{4} + j^{32} a_1^{8} + j^{34} a_1^{16} + j^{30} a_1^{32}) + x_1^{32} (j^{16} a_1 + j^{60} a_1^{2} + j^{18} a_1^{4} + j^{26} a_1^{8} + j^{29} a_1^{16} + j^{49} a_1^{32}))		$
\item{\bf	LXXX: (3)	}	$			    (x_1 (j^{41} a_1 + j^{24} a_1^{2} + j^{59} a_1^{4} + j^{8} a_1^{8} + j^{43} a_1^{16} + j^{35} a_1^{32}) + x_1^{2} (j^{34} a_1 + j^{39} a_1^{2} + j^{13} a_1^{4} + j^{35} a_1^{8} + j^{27} a_1^{16} + j^{41} a_1^{32}) + x_1^{4} (j^{18} a_1 + j^{26} a_1^{2} + j^{2} a_1^{4} + a_1^{8} + j^{2} a_1^{16} + j^{32} a_1^{32}) + x_1^{8} (j^{50} a_1 + j^{34} a_1^{2} + j^{29} a_1^{4} + j^{6} a_1^{8} + j^{40} a_1^{16} + j^{2} a_1^{32}) + x_1^{16} (j^{52} a_1 + j^{28} a_1^{2} + j^{21} a_1^{4} + j^{48} a_1^{8} + j^{4} a_1^{16} + j^{40} a_1^{32}) + x_1^{32} (j^{36} a_1 + j^{9} a_1^{2} + j^{10} a_1^{4} + j^{3} a_1^{8} + j a_1^{16} + j^{50} a_1^{32}))		$

\end{itemize}

\end{document}